\documentclass[11pt]{amsart}
\usepackage[T1]{fontenc}
\usepackage{lmodern}
\usepackage{amsmath,amssymb,mathtools,booktabs,array}
\usepackage{microtype}
\usepackage{tikz}
\usetikzlibrary{arrows.meta}
\usepackage[hidelinks]{hyperref}
\hypersetup{pdftitle={Genus monotonicity for numerical semigroups},pdfauthor={Mingchang Liu}}
\usepackage[margin=1.05in]{geometry}
\newtheorem{theorem}{Theorem}[section]
\newtheorem{lemma}[theorem]{Lemma}
\newtheorem{proposition}[theorem]{Proposition}
\newtheorem{corollary}[theorem]{Corollary}
\theoremstyle{definition}
\newtheorem{definition}[theorem]{Definition}
\theoremstyle{remark}
\newtheorem{remark}[theorem]{Remark}
\newcommand{\N}{\mathbb N_0}
\newcommand{\Fib}{\mathsf F}

\numberwithin{equation}{section}
\title{Genus monotonicity for numerical semigroups}
\author{Mingchang Liu}
\date{}
\begin{document}
\begin{abstract}
Let $n_g$ denote the number of numerical semigroups of genus $g$.
We prove that $n_{g+1}\ge n_g$ for every integer $g\ge0$, establishing
the weak monotonicity form of the Bras-Amor\'os conjecture.
The result applies to all numerical semigroups, without restrictions on
multiplicity or conductor. The computer-assisted proof begins with a
partial injection from genus $g$ to genus $g+1$. Its exceptional sources
are leaves whose largest Kunz coordinate is four or five, and an explicit family of unused
targets contains at least as many elements as this exceptional family. We establish this
comparison by combining generating-function estimates for these
semigroup populations with coefficient bounds obtained from relaxed Kunz-coordinate
conditions. The relaxed counts are evaluated by transfer recurrences on
paths and cycles. Explicit analytic estimates prove the comparison for
$g\ge837$; exact finite comparisons, combining transfer coefficients
with analytic bounds, cover $1\le g\le836$. These computations do not
require a census of all semigroups in the finite range. The stronger
Fibonacci inequality is not established by this argument.
\end{abstract}
\maketitle

\section{Introduction}

A numerical semigroup is an additive submonoid $S$ of $\N$ with finite
complement. Its multiplicity is its least positive element. Its genus is $g(S)=|\N\setminus S|$. Write
\[
 \mathcal S_g=\{S:g(S)=g\},\qquad n_g=|\mathcal S_g|.
\]
Bras-Amor\'os's computations led her to conjecture the stronger Fibonacci
inequality for this sequence~\cite{Bras}. We consider its weak monotonicity
form, recorded separately in Kaplan's survey~\cite[Conjecture~2]{Kaplan}.

\begin{theorem}\label{thm:main}
For every integer $g\ge0$, one has $n_{g+1}\ge n_g$.
\end{theorem}

The conjecture concerns all numerical semigroups, with no restriction on
multiplicity or conductor. It is weaker than the Fibonacci inequality
$n_g\ge n_{g-1}+n_{g-2}$ for $g\ge2$, which is not established here.
The distinction between all genera and sufficiently large genus is essential.
Zhai proved that
\begin{equation}\label{eq:known-asymptotic}
 n_g\sim C\varphi^g,\qquad
 \varphi=\frac{1+\sqrt5}{2},\quad C>0.
\end{equation}
Consequently $n_{g+1}/n_g\to\varphi>1$, so eventual strict increase was
already known~\cite{Zhai}. An all-genus proof requires control of the
range preceding that asymptotic conclusion. Exact values of $n_g$
have been computed through genus $80$~\cite[Section~4]{BT}.

For a semigroup of positive genus and multiplicity $m$, its Kunz word
lists, for each nonzero residue class modulo $m$, the number of gaps in
that class. The sum of its letters is the genus and the largest letter
is the depth. We use these words to compare two successive genera.

For each $g\ge1$, we construct an injection from
$\mathcal S_g\setminus\mathcal E_g$ into $\mathcal S_{g+1}$ and exhibit
at least $L(g)$ targets outside its image. Thus
\begin{equation}\label{eq:master-inequality}
 n_{g+1}-n_g\ge L(g)-|\mathcal E_g|.
\end{equation}
It therefore suffices to show that $L(g)\ge|\mathcal E_g|$. Once
there are enough unused targets, the remaining sources can be matched
to distinct targets in any fixed order. The case $g=0$ is immediate
because $n_0=n_1=1$.

Sections~\ref{cmb:injection}--\ref{cmb:binary} construct the partial
injection. Proposition~\ref{cmb:capacity} counts two disjoint supplies
of unused targets: those with binary Kunz words (letters $1$ and $2$),
and a family of nonbinary words. Their combined lower bound is $L(g)$.
It remains to count the exceptional sources. These are leaves in the
semigroup tree. From such a leaf no generator larger than its largest
gap can be deleted to give a child. Their largest Kunz coordinate is four or five.

One operation explains why two estimates of the exceptional set are useful.
Every exceptional leaf with a single five also has a four. Restoring the
Frobenius number changes the five to a four and gives a parent of genus
$g-1$ with at least two fours. Each parent has at most one leaf child
(Lemma~\ref{cmb:parent}). These parents, being nonleaves, are distinct
from the depth-four exceptional leaves. The two estimates use this observation
in different ways.
\begin{itemize}
\item For $g\ge837$, consider the sum of $z^{g(S)}$ over a family of
  semigroups at a fixed $0<z<1$. Moving a leaf to its parent increases
  its weight. A single weighted sum over depth-four words with repeated
  fours therefore bounds both families together (Lemma~\ref{cmb:compensation});
  leaves with repeated fives are counted separately. Section~\ref{mass:section}
  bounds these two sums explicitly. Proposition~\ref{mass:tail} compares
  the resulting bound on $|\mathcal E_g|$ with the sharper lower bound
  \eqref{cmb:capacity-improved} from Proposition~\ref{cmb:capacity}.
\item For $1\le g\le836$, compare counts one genus at a time. A leaf
  and its parent now contribute at different indices, so the depth-four
  terms occur at $g$ and $g-1$. Proposition~\ref{cmb:raw-envelope} bounds
  them by counting words subject to just the Kunz inequalities targeting
  their two rightmost fours, then subtracting specified families. Such
  words are called \emph{raw}; a \emph{genuine} word satisfies every Kunz
  inequality. Section~\ref{sec:computation} computes the retained coefficients,
  using Theorem~\ref{raw:tail} to bound the omitted pivot positions and a
  depth-five estimate from Section~\ref{mass:section} starting at genus $301$.
  Lemma~\ref{cmp:finite} checks the comparison for the entire finite range.
\end{itemize}

\begin{figure}[!htbp]
\centering
\begin{tikzpicture}[
  box/.style={draw,rounded corners=1pt,align=flush center,inner sep=4pt,
              font=\small,text width=5.5cm},
  wide/.style={box,text width=12.2cm},
  arr/.style={-{Stealth[length=2mm]},thin},
  reuse/.style={-{Stealth[length=2mm]},thin,densely dashed},
  lab/.style={font=\footnotesize,inner sep=1.5pt,align=flush center}]
\node[wide] (construction) at (0,0)
 {\textbf{Sections 2--4.} Partial injection $\mathcal S_g\setminus\mathcal E_g\to\mathcal S_{g+1}$
  and at least $L(g)$ unused low-depth targets
  (Propositions~\ref{cmb:partial} and~\ref{cmb:capacity})};
\node[box,minimum height=1.95cm] (raw) at (-3.55,-2.15)
 {\textbf{Section 6.} Raw two-pivot relaxation: rightmost pivots $p>300$ add at most
  $0.221715\,\mathsf F_g$ (Theorem~\ref{raw:tail}, Corollary~\ref{raw:integer-envelope})};
\node[box,minimum height=1.95cm] (mass) at (3.55,-2.15)
 {\textbf{Section 7.} Convergent masses of the depth-four and depth-five populations:\\
  $A(z_4)<2024$,\enspace $B(z_5)<113$\enspace \eqref{mass:values}};
\node[box,minimum height=2.05cm] (finite) at (-3.55,-5.45)
 {\textbf{Section 9.} Finite range $1\le g\le 836$: integer envelope of
  Proposition~\ref{cmb:raw-envelope}, exact coefficients, Lemma~\ref{cmp:finite}};
\node[box,minimum height=2.05cm] (eventual) at (3.55,-5.45)
 {\textbf{Section 8.} Eventual range $g\ge 837$: Lemma~\ref{cmb:compensation} against
  capacity \eqref{cmb:capacity-improved}, Proposition~\ref{mass:tail}};
\node[wide] (finish) at (0,-7.55)
 {$|\mathcal E_g|\le L(g)$ for every $g\ge1$; with \eqref{cmb:fundamental} this gives
  Theorem~\ref{thm:main} \\ (the case $g=0$ is checked directly)};
\draw[arr] (raw.south) -- node[lab,left=3pt] {completes $U_g$\\(from $g=604$)} (finite.north);
\draw[arr] (mass.south) -- node[lab,right=3pt] {$A$ and $B$\\masses} (eventual.north);
\draw[arr] ([xshift=5mm]mass.south west) -- node[lab,left=16pt,pos=.3]
 {$B$ gives $V_g$\\for $g\ge301$} ([xshift=-5mm]finite.north east);
\draw[arr] (construction.west) -- (-6.9,0) |- node[lab,left=2pt,pos=.25] {$L(g)$} (finite.west);
\draw[arr] (construction.east) -- (6.9,0) |- node[lab,right=2pt,pos=.25] {$L(g)$} (eventual.east);
\draw[arr] (finite.south) -- (finish.north -| finite.south);
\draw[arr] (eventual.south) -- (finish.north -| eventual.south);
\draw[reuse] (raw.east) -- node[lab,above=1pt] {\eqref{raw:pure-contraction}} (mass.west);
\end{tikzpicture}
\caption{Dependencies of the proof. Each box names its section, its output and the
 results cited there. Solid arrows are quantitative inputs. The dashed arrow is
 structural only: Section~7 reuses the carry-path decomposition and the homogeneous-pairing
 formula~\eqref{raw:pure-contraction} of Section~6, not the numerical bound of
 Theorem~\ref{raw:tail}. Here $A$ and $B$ are the generating functions of genuine
 depth-four and depth-five Kunz words with a repeated maximum (Table~\ref{tab:notation}),
 evaluated at $z_4=1247/2000$ and $z_5=633/1000$, and $U_g$, $V_g$ are the integer bounds
 entering Proposition~\ref{cmb:raw-envelope}.}
\label{fig:proof-roadmap}
\end{figure}
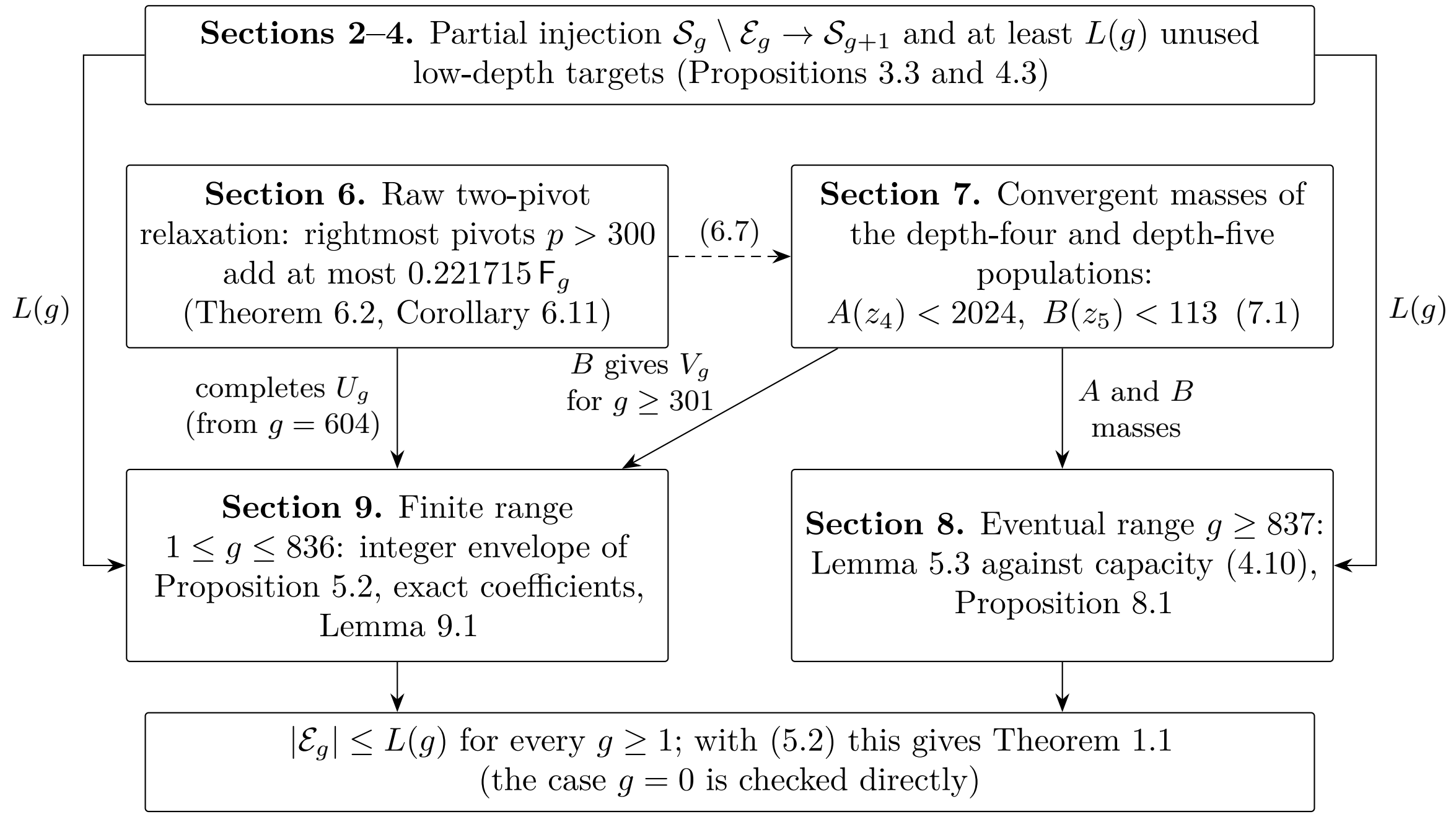

Figure~\ref{fig:proof-roadmap} summarizes these dependencies.

The proof uses established low-depth constructions. Eliahou and Fromentin
constructed two disjoint injections on gapsets of depth at most three,
which become appending $1$ and appending $2$ in Kunz coordinates~\cite[Propositions~5.2--5.3]{EF}.
Their results imply a Fibonacci lower recurrence for the low-depth count.
Kunz coordinates and their interpretation as words or compositions also
underlie the asymptotic work of Bacher and Zhu~\cite{Bacher,Zhu}.
Bacher uses maximal pivots, reflections and transfer matrices; Zhu proves
sub-Fibonacci exponential bounds for the higher-depth population using
weighted graph-homomorphism estimates. Here the estimates count the
exceptions to a particular injection, and the constants allow the finite
and analytic ranges to meet.

The all-genus monotonicity problem is distinguished explicitly from
its stronger variants in recent work of Cyrusian and Kaplan~\cite{CK}.
Bras-Amor\'os and Torra obtain refined recurrences under multiplicity
and Frobenius-number hypotheses~\cite{BT}; Casas, Madrid and Rosales
study related internal and leaf counts~\cite{CMR}.
For semigroups of ordinarization number three, Kumar's preprint combines
a shift on a large class with a rank matching into a disjoint reserve,
using separate finite and eventual capacity bounds~\cite[Sections~7--8]{Kumar}.
An earlier proof of monotonicity was proposed by Ataei Kachouei and
Rahmati~\cite{KR}. The proof of their Theorem~19 (p.~13) asserts
\[
 \binom{n_{m,g-1}}{m-1}\le n_{m,g}\qquad(n_{m,g-1}\ge m),
\]
where $n_{m,g}$ counts semigroups of multiplicity $m$ and genus $g$.
However, $n_{3,9}=n_{3,10}=4$, so the bound would require
$6=\binom42\le4$.
Indeed, multiplicity-three Kunz pairs satisfy
$a+b=g$, $2a\ge b$ and $2b+1\ge a$, giving these counts directly.

Our finite computation verifies the bounds in~\eqref{eq:master-inequality}
without a census of all numerical semigroups through genus $836$.
The population notation used throughout the proof is collected in
Table~\ref{tab:notation}. Auxiliary matrix entries and scalar bounds are
local to the subsection in which they are defined.

\begin{table}[ht]
\centering
\caption{The principal populations and counting bounds.}
\label{tab:notation}
\begin{tabular}{@{}lp{.76\linewidth}@{}}
\toprule
Notation & Meaning\\
\midrule
$\mathcal S_g$, $n_g$ & All genus-$g$ semigroups and their number\\
$\mathcal E_g$ & Exceptional leaves in the three sectors of
Lemma~\ref{cmb:residual}: depth four with at least two $4$'s and a $1$
(excluding the sole-final-$1$ case); depth five with a unique $5$, a $4$,
and a $1$; or depth five with at least two $5$'s and a $1$\\
$A_g$, $B_g$ & Genuine depth-four and depth-five words with at least two
$4$'s and at least two $5$'s, respectively; no leaf or occurrence-of-$1$
condition\\
$A(z)$, $B(z)$ & Generating series of these two genuine populations\\
$Z_{m,r,p}(z)$ & Raw depth-four partition with rightmost fours at $r<p$\\
$F_{r,p}(z)$ & Sum of $Z_{m,r,p}(z)$ over $m>p$\\
$U_g$ & Integer upper bound on the complete raw cap-four coefficient\\
$V_g$ & Integer upper bound on the genuine population $B_g$\\
$N_g$ & Exact count of $\{2,3,4\}$-words with at least two $4$'s\\
$R_g$, $L(g)$ & Exact number of unused low-depth genus-$(g+1)$ targets and
an explicit lower bound for it\\
\bottomrule
\end{tabular}
\end{table}

\section{A partial injection and its exceptional sources}
\label{cmb:injection}

This section constructs the injection on the sectors that can be handled
directly and identifies the three sectors left in $\mathcal E_g$.
Sections~\ref{cmb:binary} and~\ref{cmb:reservoir} then build disjoint reserves
of low-depth targets for the remaining matching and the exceptional sources.

We use positive Kunz coordinates~\cite{Kunz,Rosales}. If a numerical semigroup $S$ has positive
genus and multiplicity $m$, let $i+mx_i$ be its least element congruent to
$i$ modulo $m$, for $1\le i<m$. Its Kunz word is
$x=(x_1,\ldots,x_{m-1})$. The empty word represents $\mathbb Z_{\ge0}$.
The weight of a word is $|x|=\sum_i x_i$; its depth is $\max_i x_i$
for a nonempty word. The empty word has weight and depth zero.
A word is called binary if all its letters belong to $\{1,2\}$.
Throughout, coefficients with negative indices are zero.

\begin{lemma}\label{cmb:kunz}
Positive words of length $m-1$ represent numerical semigroups of
multiplicity $m$ precisely when
\begin{align}
 x_i+x_j&\ge x_{i+j} &&(i+j<m),\label{cmb:noncarry}\\
 x_i+x_j+1&\ge x_{i+j-m} &&(i+j>m),\label{cmb:carry}
\end{align}
including self-pairs. In the second inequality, adding the two residues
passes $m$, contributing the extra unit on the left; we call this a
\emph{carry}. The genus is $|x|$. If the depth is $a$ and $p$ is the
rightmost coordinate equal to $a$, the Frobenius number is $(a-1)m+p$.
\end{lemma}
\begin{proof}
In residue $i$, the gaps are $i,i+m,\ldots,i+(x_i-1)m$, proving the genus
formula. Adding the least nongaps in two nonzero residues proves necessity
of the displayed inequalities. Conversely, these inequalities ensure that
sums of these least nongaps belong to the prescribed residue tails. Adding
multiples of $m$ then proves closure for all nongaps. A sum in residue zero
already belongs to the semigroup, so no condition is needed when $i+j=m$.
The largest gap in residue $i$ is $i+(x_i-1)m$, and comparison of these
numbers proves the Frobenius formula.
\end{proof}

We set the Frobenius number of $\mathbb Z_{\ge0}$ to $-1$.
An \emph{effective generator} is a minimal generator exceeding the Frobenius
number. Deleting one gives an ordinary child of genus one greater. A
semigroup without effective generators is a \emph{leaf}. We call depth at
most three low depth and depth at least four high depth.

\begin{lemma}\label{cmb:parent}
Every positive-genus semigroup has a unique ordinary parent, obtained by
restoring its Frobenius number. Every parent has at most one leaf child.
For a high-depth word, restoring the Frobenius decreases its rightmost
maximal coordinate by one and preserves multiplicity. Consequently a
high-depth leaf has a low-depth parent exactly when its depth is four and
it has a unique coordinate equal to four.
\end{lemma}
\begin{proof}
Let $F$ be the Frobenius number. Restoring $F$ preserves closure: every sum
of $F$ and a positive nongap is larger than $F$, as is $2F$. In the enlarged
semigroup, $F$ cannot be a sum of two positive elements, for such a
decomposition would already have contradicted closure before restoration.
The new Frobenius is smaller than $F$, so $F$ is an effective generator.
Conversely, deletion of an effective generator $a$ gives a semigroup whose
Frobenius is $a$, recovering its parent uniquely.

If $a<b$ are effective generators of a parent, deleting $a$ leaves $b$
minimal, because deletion cannot introduce a decomposition. It also leaves
$b$ effective, because the child's Frobenius is $a$. Thus only deletion of
the largest effective generator can give a leaf. At high depth, $F>m$;
restoration and deletion of effective generators therefore preserve $m$.
The coordinate statement follows from Lemma~\ref{cmb:kunz}, and yields the
last assertion. At multiplicity two a word $(k)$ has effective generator
$2k+1$; in particular this boundary introduces no high-depth leaves.
\end{proof}

All sets at a fixed genus are finite. Indeed the map $s\mapsto F-s$
injects the nongaps in $[0,F]$ into its gaps, whence $F\le2g-1$ at positive
genus $g$. Ordering the finite gap sets lexicographically, for example by
their indicator vectors in $[1,2g-1]$, supplies an order whenever we use a
finite-set matching below.

Fix $g\ge1$. Partition the set $\mathcal S_g$ of genus-$g$ semigroups by
the following tests, applied in order; each test selects only sources not
selected earlier:
\begin{enumerate}
\item $\mathcal L_g$: all low-depth words;
\item $\mathcal P_g$: all high-depth nonleaves;
\item $\mathcal C_g$: all remaining high-depth leaves with a low-depth parent;
\item $\mathcal T_g$: all remaining words $(5,u)$ with every letter of $u$
  at most three;
\item $\mathcal D_g$: all remaining words of depth at least six, together
  with all remaining depth-four or depth-five words containing no $1$;
\item $\mathcal I_g$: all remaining words with a unique $5$, at a position
  at least two, and all other letters at most three;
\item $\mathcal B_g$: all remaining depth-four words whose sole $1$ is final;
\item $\mathcal E_g$: all remaining sources.
\end{enumerate}

\begin{lemma}\label{cmb:residual}
The exceptional set $\mathcal E_g$ is precisely the set of leaves containing
at least one $1$ in the following three disjoint sectors:
\begin{enumerate}
\item depth four with at least two $4$'s, excluding words whose sole $1$ is final;
\item depth five with a unique $5$ and at least one $4$;
\item depth five with at least two $5$'s.
\end{enumerate}
\end{lemma}
\begin{proof}
The first two tests leave only high-depth leaves. By Lemma~\ref{cmb:parent},
the third removes exactly the depth-four leaves with a unique $4$.
The fifth removes all depths at least six and all words without $1$.
A depth-five word with only one letter at least four has a unique $5$;
it belongs to the fourth test if this is its first letter, and otherwise
to the sixth. The seventh test makes the stated further exclusion at depth
four. Every leaf in one of the three displayed sectors misses all seven
tests, proving equality.
\end{proof}

The append-$1$ and append-$2$ maps on low-depth words are the
Eliahou--Fromentin injections~\cite[Propositions~5.2--5.3]{EF}. We recall their elementary proof
because the change of multiplicity matters here. Appending either letter
to a word of multiplicity $m$ makes the new multiplicity $m+1$. Old
noncarry inequalities remain unchanged. The new noncarry target, at index
$m$, is at most two; all its inequalities are automatic. Every new carry
left side is at least three, and every target is at most three. Thus both
append operations preserve validity.

Define the following maps:
\begin{equation}\label{cmb:maps}
\begin{array}{rcll}
 x&\longmapsto&(x,1),&x\in\mathcal L_g,\\
 x&\longmapsto&\text{the child obtained by deleting the largest effective generator},
      &x\in\mathcal P_g,\\
 x&\longmapsto&(\operatorname{par}(x),2),&x\in\mathcal C_g,\\
 (5,u)&\longmapsto&(2,u,1,3),&(5,u)\in\mathcal T_g.
\end{array}
\end{equation}
Here $\operatorname{par}(x)$ denotes the parent's Kunz word. For
$x=(w,1)\in\mathcal B_g$, replace each letter of $w$ using
\begin{equation}\label{cmb:code}
 2\longmapsto11,\qquad3\longmapsto21,\qquad4\longmapsto22,
\end{equation}
and append $2$.

Each such image is a binary word ending in $2$. The substitution is
reversible: remove the final $2$ and decode successive pairs. Thus
$|\mathcal B_g|$ distinct binary targets are already occupied before
we match $\mathcal D_g$ and $\mathcal I_g$.

\begin{lemma}\label{cmb:elementary-maps}
These maps are injective into $\mathcal S_{g+1}$ and their image sectors
are pairwise disjoint. The $\mathcal B_g$ image consists of binary words
ending in $2$.
\end{lemma}
\begin{proof}
The append-$1$ map is valid and reversible. A source in $\mathcal C_g$
has a depth-three parent of genus $g-1$; appending $2$ therefore gives a
nonbinary word of genus $g+1$. Removing this final letter determines the
parent, and Lemma~\ref{cmb:parent} determines its leaf child uniquely.
The ordinary-child map is valid, has high depth, and is inverted by
restoration of the Frobenius number.

For the fourth map, suppose the original multiplicity is $m$. The new
multiplicity is $m+2$, and the weight changes by $-5+2+1+3=1$.
Old noncarry targets have value at most three. If an addend has changed,
its new value is $2$, so its sum with any positive letter is sufficient.
The target at index $m$ is $1$. For the final target $3$ at index $m+1$,
a pair of old indices had, in the source, a carry inequality targeting the
initial $5$, so its old sum was at least four. If one of these indices were
$1$, the other would be $m$, which is not an old index. Thus these sums
are unchanged. The sole additional pair is $(1,m)$, of values $(2,1)$.
Every carry inequality is automatic at depth three. Removing the final
$1,3$ and changing the initial $2$ back to $5$ gives the inverse.

In the last map, $w$ is a word over $\{2,3,4\}$ with at least two $4$'s.
Every block preserves weight, and replacement of the final $1$ by $2$
raises genus by one. Binary words satisfy all Kunz inequalities. Remove
the final $2$ and parse consecutive pairs to recover $w$ uniquely.

The first four images are respectively low depth ending in $1$, high
depth, low depth and nonbinary ending in $2$, and low depth ending in $3$.
The code image is binary ending in $2$. These descriptions are disjoint.
\end{proof}

The direct maps have now assigned every source except
$\mathcal D_g\cup\mathcal I_g\cup\mathcal E_g$, with $\mathcal B_g$
already assigned its binary targets. The next section matches
$\mathcal D_g\cup\mathcal I_g$ into the rest of the same binary sector and
measures what remains unused.

\section{Binary allocation}
\label{cmb:binary}
Write $\Fib_0=0$, $\Fib_1=1$, and
$\Fib_{g+1}=\Fib_g+\Fib_{g-1}$ for $g\ge1$.
There are exactly $\Fib_g$ binary words of weight $g+1$ ending in $2$:
remove the final $2$ and count compositions of $g-1$ over $\{1,2\}$.
This includes $g=1$, when the remaining composition is empty.
The substitution~\eqref{cmb:code} has used $|\mathcal B_g|$ of these
targets. Can the $\Fib_g-|\mathcal B_g|$ others accommodate
$\mathcal D_g$ and $\mathcal I_g$, and how many will then remain unused?

Define the following nonnegative formal series:
\begin{align*}
 a_k(z)&=\sum_{\lceil k/2\rceil\le a\le k}z^a,&
 b_k(z)&=\sum_{\substack{1\le a,b\le k\\a+b\ge k}}z^{a+b},&
 K_k(z)&=\frac{1+a_k(z)}{1-b_k(z)},\\
 a(z)&=\frac{z^4}{1-z},&
 b(z)&=\frac{z^8(7-6z)}{(1-z)^2}.
\end{align*}
Put
\begin{align}
 M(z)&=z^6K_6K_5+z^7K_7K_6+z^8K_8K_7
       +\frac{z^9(1+a)^2}{(1-z)(1-b)^2},\label{cmb:M}\\
 H(z)&=\frac{z^4+z^5}{(1-z^2-z^3)(1-z^2-z^3-z^4-z^5)},\nonumber\\
 J(z)&=\frac{2z^{12}+4z^{14}+2z^{15}}
 {(1-2z^5-z^6)(1-z^2-z^3)(1-3z^4-2z^5-z^6)},\nonumber\\
 K(z)&=\frac{z^9}{(1-z^2-z^3)^2(1-z^2-z^3-z^4)},\nonumber\\
 W(z)&=M(z)+H(z)+J(z)+K(z),\qquad W_g=[z^g]W(z).\label{cmb:W}
\end{align}
The notation $K(z)$ and $K_k(z)$ denotes different functions.

\begin{lemma}\label{cmb:majorant}
For every $g\ge1$,
\[
 |\mathcal D_g|+|\mathcal I_g|+|\mathcal B_g|\le W_g.
\]
\end{lemma}
\begin{proof}
Let $q$ be the maximum value and $p$ its rightmost position in a word
of length $m-1$. The left block is paired by $i\mapsto p-i$: each
pair has sum at least $q$, and its letters are at most $q$. The right
block is paired by $i\mapsto m+p-i$: each pair has sum at least $q-1$,
and its letters are at most $q-1$. Here $b_k$ counts ordered pairs. A central self-pair has the corresponding
ceiling-half lower bound. For either block, summing over even and odd
lengths gives $K_q$ or $K_{q-1}$. The two block lengths determine the pivot
position, so $z^qK_qK_{q-1}$ bounds all words of depth $q$ without an
additional position factor. For fixed $q$, the factor
$z^qK_qK_{q-1}$ is Bacher's generating series for weakly admissible
compositions of maximum $q$~\cite[Proposition~9.1, equations~(10)--(12),
pp.~16--17]{Bacher}. We use this exact factor for $q=6,7,8$; the $q\ge9$
term below is a further relaxation, and the other summands of $W$ bound
different populations. For $q\ge9$,
allow arbitrary positive letters, weaken both pair thresholds to eight and
both central thresholds to four.
The central series is $a$, and the ordered pair series is
$\sum_{s\ge8}(s-1)z^s=b$. These are Zhu's single-coordinate and pair
series $\alpha$ and $\beta$ for $q_0=9$~\cite[Section~5.2]{Zhu}.
Here the two block lengths are summed separately. Summing the pivot weight over $q\ge9$ yields
the last term of~\eqref{cmb:M}. Thus $M$ bounds depth at least six.
Subtracting the composition series on $\{2,3\}$ from that on
$\{2,3,4,5\}$ gives $H$, which counts all no-$1$ words of depth four or
five. These two estimates bound $\mathcal D_g$.

For an isolated $5$ at a noninitial position, the left block is nonempty,
uses letters at most three, and has reflection-pair sum at least five.
In particular it contains no $1$. Its series is $L_0-1$, where
\[
 L_0=\frac{1+z^3}{1-2z^5-z^6}.
\]
The right block has pair threshold four and central threshold two, giving
\[
 Q_0=\frac{1+z^2+z^3}{1-3z^4-2z^5-z^6}.
\]
Since the no-$1$ words were already removed, the right block must contain
$1$. Its subpopulation without $1$ is exactly all $\{2,3\}$ words, with
series $C_0=(1-z^2-z^3)^{-1}$. Direct simplification gives
$z^5(L_0-1)(Q_0-C_0)=J$, proving the required bound for $\mathcal I_g$.

Finally, splitting a $\{2,3,4\}$ word at its first two $4$'s gives the series
\begin{equation}\label{cmb:N}
 N(z)=\frac{z^8}{(1-z^2-z^3)^2(1-z^2-z^3-z^4)}.
\end{equation}
A final $1$ supplies the factor $z$, so $K=zN$ majorizes $\mathcal B_g$.
We discarded inequalities only to enlarge each population. Summing the
three bounds proves the assertion.
\end{proof}

\begin{lemma}[Binary numerical bound]\label{cmb:binary-numerical}
For all $g\ge0$, $W_g\le13\Fib_g/25$. For $g\ge141$,
$W_g\le13\Fib_g/1000$. Moreover, with $r=809/500$,
\begin{equation}\label{cmb:binary-decay}
 \frac{W_g}{\Fib_g}<\frac{50r^2}{((63/100)r)^g}\qquad(g\ge1).
\end{equation}
\end{lemma}
\begin{proof}
The exact finite arithmetic inputs are
\begin{align*}
 25W_g&\le13\Fib_g &&(0\le g\le499),\\
 1000W_g&\le13\Fib_g &&(141\le g\le499),\\
 W(63/100)&<50,&
 \frac{50r^2}{((63/100)r)^{500}}&<\frac{13}{1000}.
\end{align*}
All geometric denominators in~\eqref{cmb:W} are positive at $63/100$.
These are finite rational checks. In particular, for each rational summand
$P(z)/Q(z)$, with $Q(0)=1$, its coefficients are calculated using
\begin{equation}\label{cmb:coefficient-recurrence}
 f_n=[z^n]P-\sum_{j=1}^{\min(n,\deg Q)}[z^j]Q\,f_{n-j}.
\end{equation}
This recurrence and the displayed rational formulas specify the exact
arithmetic checks without numerical rounding.
Since $r^2<r+1$, the two initial cases and induction in the Fibonacci
recurrence give $\Fib_g\ge r^{g-2}$ for $g\ge1$. Nonnegativity gives
$W_g(63/100)^g\le W(63/100)<50$, proving~\eqref{cmb:binary-decay}.
The base $(63/100)r$ is greater than one, so the sharper finite bound
extends to every $g\ge500$ and also implies the weaker bound there.
\end{proof}

In the fixed finite-set orders, match $\mathcal D_g\cup\mathcal I_g$
to the binary words ending in $2$ that are still available. Lemmas~\ref{cmb:majorant} and
\ref{cmb:binary-numerical} bound the combined demand
$|\mathcal D_g|+|\mathcal I_g|+|\mathcal B_g|$ by $13\Fib_g/25$, hence by
the supply $\Fib_g$.
This matching is injective and its images are disjoint from all images
in Lemma~\ref{cmb:elementary-maps}. We have proved the following.

\begin{proposition}\label{cmb:partial}
For every $g\ge1$, there is an injection
$\Phi_g:\mathcal S_g\setminus\mathcal E_g\longrightarrow\mathcal S_{g+1}$.
At least
\begin{equation}\label{cmb:binary-capacity}
 L_{\rm bin}(g)=\max\left\{\left\lceil\frac{12\Fib_g}{25}\right\rceil,
                                      \Fib_g-W_g\right\}
\end{equation}
binary targets of genus $g+1$ ending in $2$ lie outside its image.
\end{proposition}
The two entries in this maximum bound the same unused set. In particular,
they must not be added. For $g\ge141$, this gives
$L_{\rm bin}(g)\ge987\Fib_g/1000$.

This accounts for the unused part of the binary end-$2$ sector. The next
section constructs a disjoint nonbinary reserve and combines the two
capacities into $L(g)$.

\section{Witnessed unused targets}
\label{cmb:reservoir}
Among the targets already used by $\Phi_g$, only images of
$\mathcal C_g$ are nonbinary low-depth words ending in $2$. Removing
their final $2$ gives the parent of a high-depth leaf. To find unused
targets of the same form, we construct base words for which such a child
is impossible: at every position where a $3$ could become a $4$, two
letters equal to $1$ witness the Kunz inequality that the change would
violate.
Let a prefix $P$ of length $n\ge1$, ending in $3$, be followed by a binary
tail of length $\ell$. The base word has multiplicity $m=n+\ell+1$.
A prefix $3$ at position $i=n-d$ is \emph{witnessed} if two tail positions
$n+j,n+k$, with $1\le j\le k\le\ell$ and both containing $1$, satisfy
\begin{equation}\label{cmb:witness}
 j+k=\ell+1-d.
\end{equation}
We permit $j=k$. This is exactly the carry equality
$(n+j)+(n+k)=m+i$. Hence replacing that $3$ by $4$ violates its carry
inequality, whose left side is $1+1+1=3$. In semigroup terms, the
least nongap $i+3m$ in residue class $i$ is the sum of the two nongaps
$(n+j)+m$ and $(n+k)+m$. Deleting $i+3m$ would break closure under
addition.

\begin{lemma}\label{cmb:witness-unused}
If a base word is valid and every $3$ is witnessed, appending $2$ gives a
nonbinary target outside the image of $\Phi_g$, where $g$ is one greater
than the weight of the base. Distinct bases give distinct unused targets.
\end{lemma}
\begin{proof}
The base has depth three, so its Frobenius exceeds its multiplicity. Its
ordinary children preserve multiplicity, and a high-depth child would
have to increase a coordinate $3$ to $4$. Every such increase is blocked
by~\eqref{cmb:witness}, so no such child exists. Appending $2$ preserves
validity. A target in the $\mathcal C_g$ image would recover, on removing
this $2$, a base having a high-depth leaf child, a contradiction. Depth,
final letter, and the presence of $3$ exclude every other image sector.
Deleting the appended $2$ recovers the base. Its weight is $g-1$, while
the target weight is $g+1$, giving the stated genus shift.
\end{proof}

Put $u=z^2+z^3$, $s=z+z^2+z^3$, $c=z+z^2$, and $v=us$.
We use two disjoint families of prefixes. A \emph{safe-half prefix} of
length $n$ has its first $\lfloor n/2\rfloor$ letters in $\{2,3\}$,
its remaining nonfinal letters in $\{1,2,3\}$, and its last letter $3$.
This is Bacher's family in the proof of~\cite[Theorem~15.1, Section~15.1]{Bacher}.
We use it here with tail witnesses to obtain targets outside the partial
injection; the prefix family itself is not new.
An \emph{early-single-one prefix} has its unique prefix $1$ at a position
$p\ge1$, has $2$ at position $2p$, ends in $3$ at some $n\ge2p+1$,
and has all other letters in $\{2,3\}$.

Every base obtained from either family and an arbitrary binary tail is
valid. At depth three only a noncarry $1+1$ targeting $3$ can fail. In the
safe-half family two prefix $1$ positions sum beyond the prefix. In the
early-single-one family the only prefix pair of $1$'s is the self-pair at
$p$, targeting the prescribed $2$ at $2p$. A pair involving a tail position
cannot noncarry-target the prefix. These observations prove validity.
The early $1$ lies in the first half, so the two families are disjoint.
The last $3$ uniquely identifies the prefix in a base.

Writing $U=\{2,3\}$ and $V=\{1,2,3\}$, the safe prefixes have forms
$U^kV^k3$ and $U^{k+1}V^k3$, where powers denote independently chosen
letters. The early prefixes have form $U^h1U^h2U^k3$; the two $U^h$ blocks
need not be identical. Therefore their weight series are
\begin{equation}\label{cmb:prefix-series}
 C_s=\frac{z^3(1+u)}{1-v},\qquad
 C_e=\frac{z^6}{(1-u)(1-u^2)}.
\end{equation}

We now count bases with at least one unwitnessed $3$ by marking a failed
$3$. Let $t$ record the distance of the mark from the prefix end. The two
marked-prefix series are
\begin{align}
 F_s(z,t)&=C_s+
 \frac{z^6u(1+u)t}{(1-v)(1-vt)}+
 \frac{z^6t(1+st)}{(1-vt)(1-vt^2)},\label{cmb:marked-safe}\\
 F_e(z,t)&=C_e+
 \frac{z^9t}{(1-u^2)(1-u)(1-ut)}+
 \frac{z^9ut^2}{(1-ut)(1-u^2)(1-u^2t)}
 \nonumber\\
 &\quad+\frac{z^9ut^4}{(1-ut)(1-u^2t)(1-u^2t^2)}.
 \label{cmb:marked-early}
\end{align}
To verify these formulas directly, marking the terminal $3$ contributes
$C_s$ or $C_e$. In a block of length $k$ with letter-weight $a$, marking a
$3$ with $j$ letters to its right replaces its block weight by
$z^3a^{k-1}t^j$, for $0\le j<k$. Multiply further by $t$ raised to the
total length of all subsequent blocks. In the safe forms above, marking the $V^k$ blocks gives
the second term of~\eqref{cmb:marked-safe}; marking the first $U$ blocks
gives the third. In the early form, mark in order the $U^k$ block, the
second $U^h$ block, and the first $U^h$ block. These give the three
nonterminal terms of~\eqref{cmb:marked-early}, respectively. Summing each
block by the finite geometric identity
$\sum_{j=0}^{k-1}t^j=(1-t^k)/(1-t)$ yields the displayed expressions.
Thus these nonnegative series mark actual occurrences, with multiplicity
when a prefix has several $3$'s.

For a mark at distance $d$, let $Q_d(z)$ count binary tails that fail to
witness it. Then
\begin{equation}\label{cmb:failure-tail}
 Q_d(z)=\sum_{0\le\ell<d}c^\ell+c^dA_w(z),\qquad
 A_w(z)=\frac{1+z^2}{1-2z^3-z^4}.
\end{equation}
Indeed no witness is possible for $\ell<d$. For $\ell\ge d$, the first
$\ell-d$ positions are paired by reflection. A pair cannot be $11$, so its
possibilities are $12,21,22$, of weight $2z^3+z^4$; a central singleton
must be $2$. The last $d$ positions are arbitrary. The empty constrained
block when $\ell=d$ is included by the numerator of $A_w$.

For a marked series $F(z,t)=\sum_{d\ge0}F_d(z)t^d$, define
\begin{equation}\label{cmb:failure-operator}
 \mathcal U(F)=z\left\{\frac{F(z,1)-F(z,c)}{1-c}+A_wF(z,c)\right\}
             =z\sum_{d\ge0}F_d(z)Q_d(z).
\end{equation}
This is a nonnegative formal series counting failed marks. It upper bounds
bad bases, since each bad base has at least one failed mark. The extra
factor $z$ converts base weight into the source genus in
Lemma~\ref{cmb:witness-unused}. Set
\begin{align}
 T_s&=\frac{zC_s}{1-c},&T_e&=\frac{zC_e}{1-c},&
 U_s&=\mathcal U(F_s),&U_e&=\mathcal U(F_e).\nonumber
\end{align}
\begin{equation}
 L_{\rm nonbin}(g)=\max\{0,[z^g](T_s-U_s)\}
                  +\max\{0,[z^g](T_e-U_e)\}.\label{cmb:nonbinary-capacity}
\end{equation}
Separate clipping is legitimate because the two base families are disjoint.
Lemma~\ref{cmb:witness-unused} proves that this is a lower bound on the
number of unused nonbinary targets.

\begin{lemma}\label{cmb:nonbinary-numerical}
For $g\ge305$, $L_{\rm nonbin}(g)\ge(927/500)\Fib_g$.
\end{lemma}
\begin{proof}
Define nonnegative error series
\[
 E_s=\frac{(1+z)(1+z^2)^2}{1-v},\qquad
 E_e=\frac{z^4(1+z+z^2)}{(1-u)(1-u^2)},\qquad S=E_s+E_e+U_s+U_e.
\]
Direct multiplication of denominators gives
\begin{align*}
 T_s&=\frac1{1-c}-E_s,\\
 T_e&=\frac{z^4}{1-c}-E_e+
       \frac{z^9(2+2z+z^2)}{(1-u)(1-u^2)}.
\end{align*}
As $[z^j](1-z-z^2)^{-1}=\Fib_{j+1}$, the sum before clipping is at least
$\Fib_{g+1}+\Fib_{g-3}-[z^g]S$ for $g\ge4$.
All geometric denominators in $E_s,E_e,A_w,F_s(z,c),F_e(z,c)$ are
positive at $z=5/8$. Although $c(5/8)>1$, the defining nonnegative
series for $\mathcal U$ converges there: for this positive value of $c$,
$\sum_{\ell<d}c^\ell\le c^d/(c-1)$, so~\eqref{cmb:failure-operator}
is dominated by a constant times $F(z,c)$. We are not expanding
$(1-c)^{-1}$ as a convergent geometric series at this point. Its apparent
singularity at $c=1$ in~\eqref{cmb:failure-operator} is removable: the
quotient comes from the finite sums $\sum_{0\le\ell<d}c^\ell$.

The exact rational inputs, computed by~\eqref{cmb:coefficient-recurrence}, are
\begin{align*}
 S(5/8)&<78,\\
 78(809/500)^2(800/809)^{1500}&<1/10000,\\
 10000(\Fib_{g+1}+\Fib_{g-3})&\ge18541\Fib_g &&(g=16,17),\\
 500L_{\rm nonbin}(g)&\ge927\Fib_g &&(305\le g\le1499).
\end{align*}
The difference in the third line satisfies the Fibonacci recurrence, so
the inequality holds for every $g\ge16$. Nonnegative coefficients and
$\Fib_g\ge(809/500)^{g-2}$ give
\[
 \frac{[z^g]S}{\Fib_g}<78(809/500)^2(800/809)^g<\frac1{10000}
 \qquad(g\ge1500).
\]
Subtracting this error gives $18540\Fib_g/10000=927\Fib_g/500$ for
$g\ge1500$. The finite interval in the last line completes the proof.
\end{proof}

\begin{proposition}\label{cmb:capacity}
Define the exact integer lower bound
\begin{equation}\label{cmb:L}
 L(g)=L_{\rm bin}(g)+L_{\rm nonbin}(g).
\end{equation}
For every $g\ge1$, at least $L(g)$ low-depth targets of genus $g+1$ are outside the image of
$\Phi_g$. For $g\ge305$,
\begin{align}
 L(g)&\ge\frac{2841}{1000}\Fib_g,\label{cmb:capacity-constant}\\
 \frac{L(g)}{\Fib_g}&\ge1-
 \frac{50(809/500)^2}{((63/100)(809/500))^g}+\frac{927}{500}.
 \label{cmb:capacity-improved}
\end{align}
\end{proposition}
\begin{proof}
The binary and nonbinary unused sets are disjoint. Thus their bounds add.
The first estimate combines $987/1000$ and $927/500$; the second combines
$L_{\rm bin}(g)\ge\Fib_g-W_g$ with~\eqref{cmb:binary-decay} and
Lemma~\ref{cmb:nonbinary-numerical}.
\end{proof}

Proposition~\ref{cmb:capacity} supplies the explicit reserve needed in the
master inequality. The next section distinguishes its lower bound $L(g)$
from the exact unused low-depth count $R_g$ and isolates the remaining task
of bounding $|\mathcal E_g|$.

\section{The remaining counting problem}
\label{cmb:envelope}
The partial injection and the target reserve are now complete. What remains
is to compare the three exceptional sectors with that reserve, first by an
integer coefficient envelope and later by a positive-real mass comparison.
Let $t_g$ count low-depth genus-$g$ words, and write
$c_g=|\mathcal C_g|$, $f_g=|\mathcal T_g|$,
$d_g=|\mathcal D_g|$, $i_g=|\mathcal I_g|$, $b_g=|\mathcal B_g|$.
Only the ordinary-child map uses high-depth targets. The exact number of
unused low-depth targets is therefore
\begin{equation}\label{cmb:exact-capacity}
 R_g=t_{g+1}-t_g-c_g-f_g-d_g-i_g-b_g\ge L(g).
\end{equation}
Counting the $n_g-|\mathcal E_g|$ matched targets and these unused targets
proves the fundamental inequality
\begin{equation}\label{cmb:fundamental}
 n_{g+1}-n_g\ge R_g-|\mathcal E_g|\ge L(g)-|\mathcal E_g|.
\end{equation}
There may also be unused high-depth targets, which only strengthen it.

Let $A_g$ count all genuine depth-four words with at least two $4$'s, and
let $B_g$ count all genuine depth-five words with at least two $5$'s.
No leaf or $1$ condition is imposed in these two definitions.

\begin{lemma}\label{cmb:parent-reduction}
The unique-$5$ sector of $\mathcal E_g$ injects into nonleaf depth-four
words of genus $g-1$ with at least two $4$'s and at least one $1$.
Consequently $|\mathcal E_g|\le A_g+A_{g-1}+B_g$.
\end{lemma}
\begin{proof}
Restore the Frobenius, changing the unique $5$ to $4$. The result is its
ordinary parent, has one new and at least one preexisting $4$, and retains
its $1$'s. It is a nonleaf because the source is its child. Two such sources
cannot have the same parent by Lemma~\ref{cmb:parent}. The other two
sectors of Lemma~\ref{cmb:residual} are directly counted by $A_g$ and $B_g$.
\end{proof}

For $q\in\{4,5\}$ define the \emph{raw cap-$q$ population} as follows.
Choose a word with at least two $q$'s and denote their two rightmost
positions by $r<p$. Coordinates before $r$ belong to $\{1,\ldots,q\}$;
all nonpivot coordinates after $r$ belong to $\{1,\ldots,q-1\}$.
Retain only the Kunz inequalities targeting $r$ or $p$, including
self-pairs with one of those two targets. Every genuine word counted by $A_g$ or $B_g$ belongs to the
corresponding raw population. The two rightmost maxima are unique, so
summing over their positions introduces no multiple counting. The complete
raw coefficient sums over all multiplicities and both pivot positions.

\begin{proposition}\label{cmb:raw-envelope}
Let $U_g$ be any integer upper bound on the complete raw cap-four
coefficient of weight $g$, and let $V_g$ be any integer upper bound on
$B_g$. With $N_g=[z^g]N(z)$ from~\eqref{cmb:N}, the quantity
\begin{equation}\label{cmb:Eupper}
 E_{\rm upper}(g)=U_g-N_g-N_{g-3}-N_{g-4}
                         +U_{g-1}-N_{g-1}+V_g
\end{equation}
is an upper bound on $|\mathcal E_g|$.
\end{proposition}
\begin{proof}
Raw cap-four words without $1$ are exactly the $\{2,3,4\}$ words with
at least two $4$'s, and every such word is genuine: noncarry sums are at
least four and carry sums at least five. Their exact series is $N(z)$.
For a word whose sole $1$ is final, all inequalities are automatic except
its self-pair, which is a carry pair targeting the penultimate coordinate
and forces that coordinate to be at most three. If the penultimate
coordinate were $4$, it would be the rightmost maximum, so the raw
population retains precisely this obstruction. Thus raw and genuine
populations of this description coincide and consist of
$(v,2,1)$ and $(v,3,1)$ with $v$ a $\{2,3,4\}$ word containing at least
two $4$'s. Their series is $(z^3+z^4)N(z)$, disjoint from the no-$1$
population. Subtracting these exact counts from $U_g$ bounds the
first residual sector.

Lemma~\ref{cmb:parent-reduction} bounds the second sector by the
raw cap-four count of weight $g-1$ containing $1$, hence by
$U_{g-1}-N_{g-1}$. Its image need not be a residual leaf and may have its
sole $1$ final, so there is no further subtraction from this shifted term.
Finally $V_g$ bounds the third sector. Adding these bounds proves the
claim. The exact exclusions are from the whole raw population; the way an
upper bound $U_g$ is assembled does not affect this subtraction.
\end{proof}

For example, $(5,4,2,1)$ is a leaf with parent $(4,4,2,1)$, whose sole $1$
is final. The latter has effective generator $21$, and its deletion gives
the former. Thus the absence of a final-$1$ exclusion in the shifted term
of~\eqref{cmb:Eupper} is essential.

We will also use the following real-variable improvement. Let $D(z)$ count
the first two residual sectors together, with their actual genus weights,
and put $A(z)=\sum_gA_gz^g$, $B(z)=\sum_gB_gz^g$.
\begin{lemma}\label{cmb:compensation}
At every $0<z\le1$ for which $A(z)$ converges, $D(z)\le A(z)$.
In particular, whenever $0<z_4,z_5\le1$ and both masses are finite,
\begin{equation}\label{cmb:mass-extraction}
 |\mathcal E_g|\le A(z_4)z_4^{-g}+B(z_5)z_5^{-g}.
\end{equation}
\end{lemma}
\begin{proof}
Map each depth-four residual leaf to itself, and each unique-$5$ residual
leaf to its parent as in Lemma~\ref{cmb:parent-reduction}. The former images
are leaves and the latter are nonleaves, so they are disjoint. This gives
an injection into the population counted by $A$. The first map preserves
genus, while the second lowers it by one. For $0<z\le1$ each source weight
is therefore at most its image weight. Nonnegative summation proves
$D(z)\le A(z)$. Extract a single nonnegative coefficient at $z_4$, and bound
the remaining multiple-$5$ sector using $B$ at $z_5$, to obtain the result.
This is a positive-real comparison, not a coefficientwise assertion
$[z^g]D\le A_g$.
\end{proof}

The two branches are now reduced to their quantitative inputs.
Section~\ref{raw:section} bounds the omitted raw cap-four pivots needed to
make $U_g$ explicit; Section~\ref{mass:section} then bounds the genuine
depth-four and depth-five masses used in the finite and eventual comparisons.

\section{A coefficient bound for the raw two-pivot tail}\label{raw:section}

The finite comparison in Section~\ref{sec:computation} counts depth-four
words by the positions $r<p$ of their two rightmost fours. Its table stops
at position $p=300$. A word whose last four occurs farther to the right
is absent from that table. Theorem~\ref{raw:tail} bounds the total
contribution of all such omitted positions at each genus, even if one
imposes only the Kunz inequalities targeting $r$ or $p$. Thus the computed
prefix plus a single error bound gives the upper bound $U_g$ in
Proposition~\ref{cmb:raw-envelope}. The error allowance in $U_g$ starts
at genus $604$ and in $U_{g-1}$ at $605$. Proposition~\ref{mass:tail}
uses a different bound for $g\ge837$.

Here is how the coefficient estimate is obtained. Fixing $r,p$ turns each
retained inequality into an edge constraint. Summing coordinate weights
along the resulting paths gives a recurrence as the multiplicity varies
(Lemma~\ref{raw:rationality}). A recurrence for each pivot pair alone
would give no bound on their sum. Proposition~\ref{raw:cancellation}
isolates the Fibonacci factor common to these generating functions;
Lemmas~\ref{raw:wiener} and~\ref{raw:signed} control the remaining
coefficients in absolute value. The decay in the pivot positions proved
in Proposition~\ref{raw:pointwise} then makes the sum over $p>300$
small enough.

Both Bacher and Zhu turn constraints on Kunz coordinates into weighted
graphs. For reflection paths the closest precedents are Bacher's two
rightmost maximal fives and his noncarry paths for cap-four left factors
\cite[Propositions~11.3 and~12.4]{Bacher}. Zhu groups genuine Kunz words
by their dense depth, standardizes them, and bounds the resulting weighted
graph homomorphism counts \cite[Section~5]{Zhu}. The dense depth selects
a frequently occurring letter, which can be smaller than the maximum.
Here the two positions are the actual rightmost maxima, and we drop every
inequality with a different target. This admits words that need not be
semigroups. The estimate below bounds their coefficients uniformly in genus
after summing all omitted positions and multiplicities.

\begin{definition}[Raw two-pivot polynomial]\label{raw:definition}
For integers $1\le r<p<m$ and $a\in\{4,5\}$, let
$Z^{(a)}_{m,r,p}(z)$ count words $(x_1,\ldots,x_{m-1})$ with
\[
 x_r=x_p=a,\qquad
 1\le x_i\le a\ (i<r),\qquad
 1\le x_i\le a-1\ (i>r,\ i\ne p),
\]
weighted by $z^{\sum_i x_i}$.  Retain exactly the inequalities
\[
 x_i+x_j\ge a\quad(i+j\in\{r,p\}),\qquad
 x_i+x_j\ge a-1\quad(i+j\in\{m+r,m+p\}),
\]
including the cases $i=j$. No inequality with another target is imposed.
A \emph{partition} here means the sum of the weights over the admissible
coordinate assignments. Write $Z_{m,r,p}=Z^{(4)}_{m,r,p}$ and
\[
 F_{r,p}(z)=\sum_{m>p}Z_{m,r,p}(z).
\]
\end{definition}

The pivots are uniquely the two rightmost occurrences of the maximum.
Consequently summing over $m,r,p$ counts each raw word once.  A genuine
depth-$a$ Kunz word with at least two maxima belongs to this relaxation.

For a genuine depth-four word, the two largest gaps are $3m+p$ and
$3m+r$. Indeed, the gaps in residue class $i$ end at $i+(x_i-1)m$.
Only a four gives a gap above $3m$, and $p,r$ are the two rightmost fours.
Closure forbids two nongaps from adding to either of these gaps. Reducing
the two sums modulo $m$ gives precisely the reflection pairs in the
definition. For a raw word the same residues specify which inequalities
we test, even when closure fails elsewhere.

\begin{theorem}[Raw omitted-pivot estimate]\label{raw:tail}
For every integer $g\ge0$,
\begin{equation}\label{raw:main}
 \sum_{p>300}\sum_{r=1}^{p-1}[z^g]F_{r,p}(z)
 \le \frac{221715}{1000000}\,\Fib_g.
\end{equation}
The left side vanishes for $g<604$.  No leaf condition or occurrence of a
coordinate equal to one is required.
\end{theorem}

\subsection{The reflection graph and its support}

Remove the two fixed pivots and use the remaining coordinates as vertices,
with the alphabets in
Definition~\ref{raw:definition} and vertex weight $z^{x_i}$ at coordinate
$i$.  For a target $t\in\{r,p\}$, the retained pairs reflect $i$ to $t-i$
in the noncarry range and to $m+t-i$ in the carry range; this is why we call
the result a \emph{reflection graph}.
A noncarry pair gives a \emph{strong} edge,
of threshold $a$, and a carry pair gives a \emph{carry} edge, of threshold
$a-1$.  A self-pair is a unary restriction, represented as a loop rather
than an ordinary edge.  Each vertex has at most one reflection incidence
for each target.  Thus the ordinary graph has maximum degree two and its
components are paths, isolated vertices, and cycles.  A cycle alternates
the two targets modulo $m$, has even length, and has length at least four:
distinct targets cannot give parallel edges.  Constraints incident to a
fixed pivot are automatic.

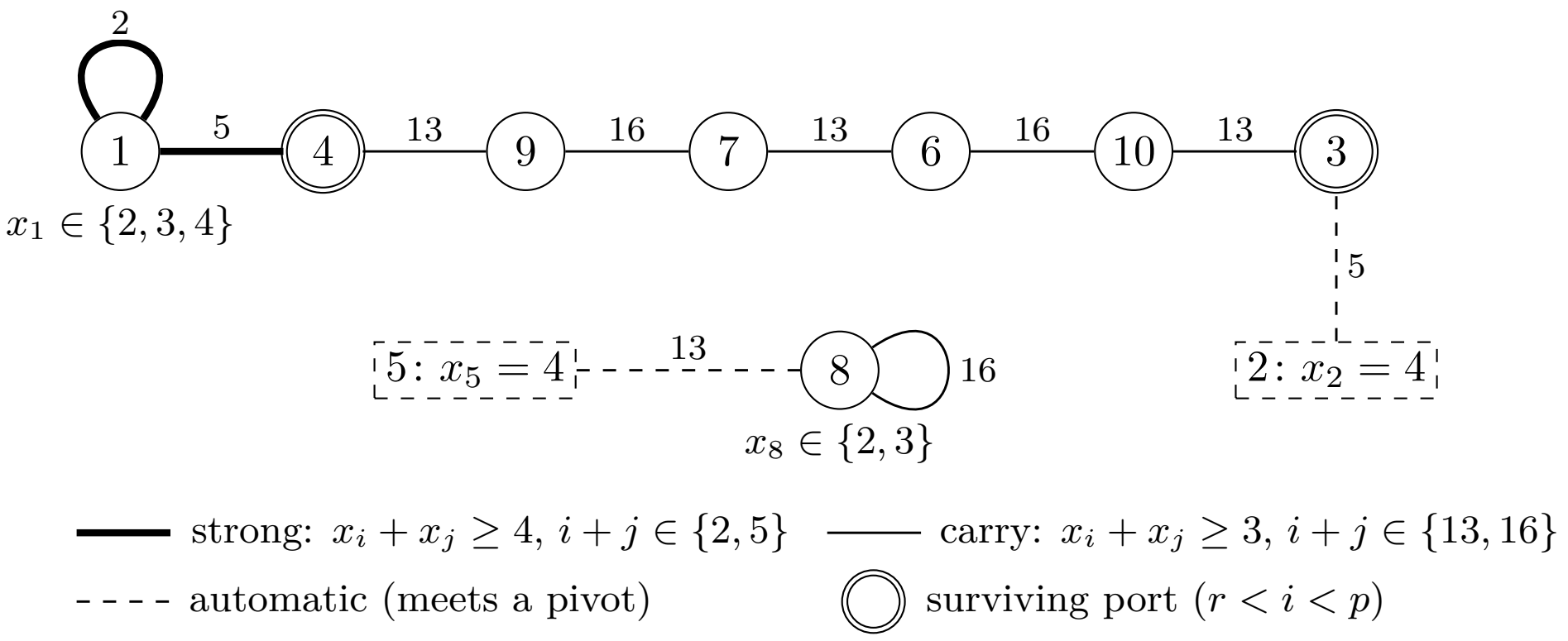
\begin{figure}[htbp]
\centering
\begin{tikzpicture}[x=1cm,y=1cm,
  v/.style={circle,draw,minimum size=6.2mm,inner sep=0pt,font=\small},
  port/.style={v,double,double distance=.8pt},
  piv/.style={rectangle,draw,dashed,inner sep=2.5pt,font=\small},
  strong/.style={line width=1.6pt},
  carry/.style={line width=.55pt},
  auto/.style={dashed,line width=.55pt},
  el/.style={font=\scriptsize,inner sep=1.2pt,fill=white},
  note/.style={font=\footnotesize,align=left,inner sep=1pt}]
\node[v] (n1) at (0.000,0.000) {$1$};
\node[port] (n4) at (1.620,0.000) {$4$};
\node[v] (n9) at (3.240,0.000) {$9$};
\node[v] (n7) at (4.860,0.000) {$7$};
\node[v] (n6) at (6.480,0.000) {$6$};
\node[v] (n10) at (8.100,0.000) {$10$};
\node[port] (n3) at (9.720,0.000) {$3$};
\node[v] (n8) at (5.751,-1.750) {$8$};
\node[piv] (p5) at (2.835,-1.750) {$5\colon x_5=4$};
\node[piv] (p2) at (9.720,-1.750) {$2\colon x_2=4$};
\draw[strong] (n1) -- node[el,above=1pt] {$5$} (n4);
\draw[carry] (n4) -- node[el,above=1pt] {$13$} (n9);
\draw[carry] (n9) -- node[el,above=1pt] {$16$} (n7);
\draw[carry] (n7) -- node[el,above=1pt] {$13$} (n6);
\draw[carry] (n6) -- node[el,above=1pt] {$16$} (n10);
\draw[carry] (n10) -- node[el,above=1pt] {$13$} (n3);
\draw[strong] (n1) to[out=125,in=55,looseness=7] node[el,above] {$2$} (n1);
\draw[carry] (n8) to[out=35,in=-35,looseness=7] node[el,right=1pt] {$16$} (n8);
\draw[auto] (p2) -- node[el,right=1pt] {$5$} (n3);
\draw[auto] (p5) -- node[el,above=1pt] {$13$} (n8);
\node[note,below=2pt] at (n1.south) {$x_1\in\{2,3,4\}$};
\node[note,below=2pt] at (n8.south) {$x_8\in\{2,3\}$};
\draw[strong] (-0.35,-3.05) -- ++(0.75,0) node[note,right=3pt] {strong: $x_i+x_j\ge4$, $i+j\in\{2,5\}$};
\draw[carry] (5.65,-3.05) -- ++(0.75,0) node[note,right=3pt] {carry: $x_i+x_j\ge3$, $i+j\in\{13,16\}$};
\draw[auto] (-0.35,-3.60) -- ++(0.75,0) node[note,right=3pt] {automatic (meets a pivot)};
\node[port,minimum size=4.6mm] at (6.02,-3.60) {};
\node[note,right] at (6.40,-3.60) {surviving port ($r<i<p$)};
\end{tikzpicture}
\caption{The retained pairs for $a=4$, $m=11$, $r=2$, $p=5$, drawn as one graph on the nonpivot coordinates; edge labels are $i+j$. Thick edges are strong (noncarry) pairs and thin edges are carry pairs; dashed incidences meet a fixed pivot, are automatic and are deleted with the pivots. The surviving ordinary path $1\!-\!4\!-\!9\!-\!7\!-\!6\!-\!10\!-\!3$ has one strong edge and five carry edges whose index sums alternate between $m+r=13$ and $m+p=16$. The loops impose $x_1\ge2$ and $x_8\ge2$; all other nonpivot letters lie in $\{1,2,3\}$. Deleting the pivot $5$ leaves $8$ as a zero-ordinary-edge free-to-loop component of weight $z^2+z^3$. The temporary carry ports are $3,4,5$; only $3,4$ survive. Multiply the two component polynomials by the pivot factor $z^8$. The multiplicity-progression argument is Lemma~\ref{raw:rationality}.}
\label{fig:raw-constraint-graph}
\end{figure}

This description also specifies the exact polynomial: sum the products
of vertex weights over allowed assignments along each path; for a cycle,
fix one starting value and close the last edge back to that value.
Multiply the component polynomials and the pivot factor $z^{2a}$.
Loops impose their unary restrictions at the corresponding endpoints.
An isolated vertex is a path of zero ordinary-edge length.
Figure~\ref{fig:raw-constraint-graph} shows how the two types of retained
pairs fit together in one small example.

\begin{lemma}[Finite support]\label{raw:support}
Every word counted by $Z^{(a)}_{m,r,p}$ has genus satisfying
\begin{equation}\label{raw:genus-support}
 2g\ge(a-1)m+p+1.
\end{equation}
In particular, for cap four, $g\ge2p+2$ and $g\ge m+5$.
\end{lemma}
\begin{proof}
The reflection targeting $p$ pairs the $p-1$ coordinates to its left
with sum at least $a$ and the $m-p-1$ coordinates to its right with sum
at least $a-1$.  A central fixed point satisfies twice its coordinate
at least the same threshold.  Adding these inequalities, with twice the
pivot weight, gives
\[
 2g\ge a(p-1)+(a-1)(m-p-1)+2a=(a-1)m+p+1.
\]
For $a=4$ use $m\ge p+1$.  Independently, the two pivot entries contribute
eight and the other $m-3$ coordinates contribute at least one each,
giving $g\ge m+5$.
\end{proof}

All series and coefficient sums below are therefore well defined formally.
The support assertion in Theorem~\ref{raw:tail} follows already from
$p\ge301$.

\subsection{Fixed-pivot recurrence and reduction of the homogeneous denominator}\label{raw:topology}

We first fix the pivot positions and let the multiplicity vary.  The
noncarry constraints stay in a fixed finite graph; the carry components
lengthen in a controlled way.  The resulting recurrence initially has a
factor for every choice of an eigenvalue along the carry components.
The two homogeneous contributions, which use the same eigenvalue on every
component, can be summed more efficiently across multiplicities. This
regrouping removes a separation-dependent geometric factor from the
denominator; the reduced denominator admits the uniform bounds needed in
Lemmas~\ref{raw:wiener} and~\ref{raw:signed}. Neither homogeneous
contribution is discarded; the mixed factors remain.

\paragraph{\textit{Carry paths and their attachment coordinates.}}

Fix $r<p$, and put $d=p-r$.  For $m>2p$, temporarily include the pivot
$p$ as a vertex of the carry graph, still with its value fixed at $4$.
This temporary inclusion describes incidences, not an additional variable
with alphabet $\{1,2,3\}$.  Set $x=i-r$ and $N=m-r$.
The carry vertices are $1,\ldots,N-1$, and their two reflections are
\begin{equation}\label{raw:reflections}
 \mathcal F(x)=N-x,\qquad \mathcal G(x)=N+d-x.
\end{equation}
The first is defined throughout this interval.  The second stays in the
interval exactly when $x>d$.  Thus the vertices $1,\ldots,d$ have one
carry incidence each and form the \emph{ports}.  The pivot $p$ is port
$d$; the surviving ports correspond to coordinates $r+1,\ldots,p-1$.
Every carry vertex other than the fixed pivot is a variable with alphabet
$\{1,2,3\}$; this includes all surviving ports. Port $d$ is the pivot,
fixed at $4$. A port is defined by its carry incidence, not by its degree
in the full reflection graph.

Here is the complete path decomposition, including the loops.  For a
representative $a\in\{1,\ldots,d\}$, list its residue class in the interval as
$a_k=a+kd$, $0\le k<n_a$.  Let $b\in\{1,\ldots,d\}$ be determined by
$a+b\equiv N\pmod d$.  The map $\mathcal F$ reverses the two lists and
bijects them, so $n_a=n_b=n$.  Since $a+(n-1)d+b=N$, their index formulas are
\[
 \mathcal F(a_k)=b_{n-1-k},\qquad
 \mathcal G(a_k)=b_{n-k}\quad(1\le k<n).
\]
For $a\ne b$, alternating these formulas gives the path
\[
 a_0,b_{n-1},a_1,b_{n-2},\ldots,a_{n-1},b_0.
\]
It contains all $2n$ vertices in the two classes, has ports $a,b$, and
has $2n-1$ ordinary edges.  For $a=b$, the same alternation visits the
indices $0,n-1,1,n-2,\ldots$ until all $n$ vertices have appeared.  It
has one port, $n-1$ ordinary edges, and a loop at its other endpoint.
That loop belongs to $\mathcal F$ when $n$ is odd and to $\mathcal G$
when $n$ is even.  Both loops require the coordinate value to be at
least two.  The orbits of $a\mapsto N-a\pmod d$ exhaust the residue
classes, so this list omits no carry component.

Replacing $m$ by $m+d$ leaves the port involution unchanged and adds one
vertex to every residue class.  Each paired-class path gains two
ordinary edges, and each fixed-class path gains one.  If $h_c$ is the
number of paired-class paths and $f$ the number of fixed classes, then
\[
 2h_c+f=d,
\]
so the total carry-edge increment is exactly $d$.  A change of loop
target does not change its unary restriction.  Deleting port $d$ removes
an endpoint and its incident ordinary edge, never an internal vertex.
It creates one free endpoint.  The strong graph, its loops and the
alphabets on all surviving ports are independent of $m$.

Write
\[
 V_0=\{1,\ldots,p-1\}\setminus\{r\},\qquad
 P_0=\{r+1,\ldots,p-1\}.
\]
The \emph{fixed noncarry graph} (also called the strong core) has vertex
set $V_0$, all surviving strong edges and loops, and the original
alphabets.  It includes isolated vertices.  Its attachment coordinates
are exactly the surviving ports $P_0$.  Reconstructing the full graph
means multiplying its constraint and weight factors by the carry-path
kernels, then summing over all allowed coordinate values, including the
values at $P_0$.

\paragraph{\textit{Weights and zero-edge boundaries.}}
\begin{samepage}
A carry edge between two coordinates with values in $\{1,2,3\}$
allows every pair except $(1,1)$. Giving each endpoint half its exponent
produces the symmetric transfer matrix
\[
 C(z)=\begin{pmatrix}
 0&z^{3/2}&z^2\\
 z^{3/2}&z^2&z^{5/2}\\
 z^2&z^{5/2}&z^3
 \end{pmatrix},\qquad
 C(z)_{ij}=z^{(i+j)/2}{\bf1}_{i+j\ge3}.
\]
For example, on the path $4,9,7,6,10,3$ in
Figure~\ref{fig:raw-constraint-graph}, fix the endpoint values
$x_4=u$ and $x_3=v$. The entry $(C^5)_{uv}$ sums over the four
interior values and supplies half the weight at each endpoint. The
fixed noncarry graph supplies the other half, $z^{(u+v)/2}$. Multiply
by its other vertex weights and constraint indicators, then sum over
the permitted core values, including $u,v$. The separate loop component
and the fixed pivots contribute their own factors.
\end{samepage}

The matrix $C$ has rank two. Its nonzero eigenvalues $\lambda,\mu$ satisfy
\begin{equation}\label{raw:eigenvalues}
 \lambda+\mu=b=z^2+z^3,\qquad
 \lambda\mu=c=-z^3-z^4.
\end{equation}
This nonzero characteristic polynomial also occurs in Bacher's
matrix~(27) and Remark~12.11 (p.~33)~\cite{Bacher}.  That local transfer
correspondence does not establish the denominator reduction proved below.
In a product of symmetric edge matrices, each internal coordinate receives
its full weight $z^i$, while an endpoint receives only $z^{i/2}$.
An endpoint vector or a core weight supplies the other half.
For a positive ordinary-edge length $\ell$, its path kernel is
$C^\ell=\lambda^\ell P_\lambda+\mu^\ell P_\mu$.
There is one possible zero-edge case to check.  Since $m>2p$, every
residue class has at least two vertices.  Consequently a paired-class
path still has positive length after deleting $p$, whereas a fixed-class
path with two vertices becomes one free, loop-constrained vertex.
The free and loop boundary vectors for this component are
\[
 h=(z^{1/2},z,z^{3/2})^T,\qquad h_+=(0,z,z^{3/2})^T.
\]
They both belong to the image of $C$, which is the span of
$(1,0,0)^T$ and $(0,z,z^{3/2})^T$.  Therefore
\[
 h^T(P_\lambda+P_\mu)h_+=h^Th_+.
\]
The zero-edge component is thus also exactly the sum of its two
nonzero modes.  It introduces no extra transient term into the
recurrence.  This is a boundary-weighted identity; it does not assert
$C^0=P_\lambda+P_\mu$ as matrices.

\paragraph{\textit{The fixed-pivot recurrence.}}
\begin{lemma}[Fixed-pivot rationality]\label{raw:rationality}
On each progression $m=m_0+kd$ with $m_0>2p$, the polynomials
$Z_{m,r,p}$ satisfy the scalar recurrence with characteristic polynomial
\[
 P_d(T,z)=\prod_{j=0}^d(T-\lambda^{d-j}\mu^j)\in\mathbb Z[z,T].
\]
Consequently $F_{r,p}$ is a rational function.
\end{lemma}
\begin{proof}
On this progression the port pairing, endpoint restrictions and fixed
noncarry graph are unchanged.  Choose one spectral summand on each carry
path, multiply its endpoint entries by the remaining core weights and
constraint indicators, and sum over every permitted core coordinate
value, including all surviving port values.  Sum free and loop endpoint
values using their boundary vectors.  The factors independent of the
progression index give a scalar amplitude for each such choice.
A paired-class component
assigned the $\lambda$ mode contributes a factor $\lambda^{2k}$ under
the shift $m_0\mapsto m_0+kd$; a fixed-class component contributes
$\lambda^k$.  The same statements hold with $\mu$ in place of $\lambda$.
If $j$ is the sum of the increments assigned to $\mu$, the remaining
increments sum to $d-j$.  Collecting the finitely many mode assignments
therefore gives
\[
 Z_{m_0+kd,r,p}=\sum_{j=0}^d A_j(z)(\lambda^{d-j}\mu^j)^k.
\]
This identity holds over the corresponding algebraic function field,
including at $k=0$ by the zero-edge calculation above; its amplitudes
need not be positive.  Applying $\prod_{j=0}^d(T-\lambda^{d-j}\mu^j)$
to this exponential sum proves the recurrence.
For integrality of the annihilator, put
$S_0=2$, $S_1=b$, and $S_n=bS_{n-1}-cS_{n-2}$.
Pair the roots indexed by $j,d-j$ to obtain
\[
 T^2-c^jS_{d-2j}T+c^d\qquad(0\le j<d/2),
\]
with the additional factor $T-c^{d/2}$ when $d$ is even.

For an explicit rational construction define
\[
 Q_d(y,z)=\prod_{j=0}^d(1-\lambda^{d-j}\mu^j y)
          =\sum_{j=0}^{d+1}q_j(z)y^j.
\]
For $1\le s\le d$ set $m_s=2p+s$ and
\[
 N_s(y,z)=\sum_{k=0}^d\sum_{j=0}^k
              q_j(z)Z_{m_s+(k-j)d,r,p}(z)y^k.
\]
The recurrence yields
\begin{equation}\label{raw:initial-data}
 F_{r,p}(z)=\sum_{m=p+1}^{2p}Z_{m,r,p}(z)
              +\frac{\sum_{s=1}^dN_s(1,z)}{Q_d(1,z)}.
\end{equation}
The substitution $y=1$ is justified coefficientwise by $g\ge m+5$,
and $Q_d(1,0)=1$.  In particular this is a structural recurrence,
not a recurrence fitted to finitely many values.
\end{proof}

\paragraph{\textit{Homogeneous contributions and denominator reduction.}}
The terms with every carry component assigned $\lambda$, or every one
assigned $\mu$, will be called \emph{homogeneous}.  All other assignments
are genuinely mixed.  The next proposition sums the homogeneous terms
across all multiplicities before comparing their denominator with the
one furnished by the recurrence.

\begin{proposition}[Cancellation and numerator degree]\label{raw:cancellation}
Let $\tau=1$ when $d$ is odd and $\tau=2$ when $d$ is even.  Then
\begin{equation}\label{raw:correction}
 H_{r,p}:=(1-z-z^2)F_{r,p}=\frac{A_{r,p}}{B_d},\qquad
 B_d=K_\tau R_d,
\end{equation}
where $A_{r,p}\in\mathbb Q[z]$, $\deg A_{r,p}\le3d^2+7p+6$, and
\begin{align*}
 K_1&=1+z+z^2,\\
 K_2&=1+z+2z^2+z^3-z^5-z^6,\\
 R_d&=\prod_{j=1}^{d-1}(1-\lambda^{d-j}\mu^j)\in\mathbb Z[z].
\end{align*}
\end{proposition}
\begin{proof}
We first separate the multiplicity dependence of a homogeneous contribution
from the way in which its ports are paired.  Fix $\theta\in\{\lambda,\mu\}$,
choose a nonzero eigenvector $v_\theta$ of the symmetric carry matrix,
and put $t_\theta=v_\theta^Tv_\theta$.  Over an algebraic function field
these choices can be made with $t_\theta\ne0$, and
\[
 P_\theta=\frac{v_\theta v_\theta^T}{t_\theta}.
\]
Let $\Omega_0$ consist of assignments to $V_0$ in the original alphabets
that satisfy all strong edges and strong loops.  The contraction of the
fixed noncarry graph is the explicit finite sum
\[
 \mathcal C_{r,p}(v_\theta)
 =z^8\sum_{x\in\Omega_0}
       \left(\prod_{i\in V_0\setminus P_0}z^{x_i}\right)
       \left(\prod_{i\in P_0}z^{x_i/2}(v_\theta)_{x_i}\right).
\]
In particular, the surviving port values are summed, not held fixed.
The factor $z^8$ accounts for the pivots.  At each surviving port,
the other half of its weight is already in the adjacent symmetric
carry kernel.  At every other core vertex the full weight is used.
Thus this is one fixed scalar, independent of $m$ and of the port pairing.

There are $d-1$ surviving ports, one free carry endpoint created by
deleting $p$, and $f$ carry-loop endpoints.  There are
$(d+f)/2$ carry components, each contributing one rank-one
projector.  Factoring its two endpoint entries shows that the entire
homogeneous $\theta$ contribution is exactly
\begin{equation}\label{raw:pure-contraction}
 \theta^{e_C(m)}\mathcal C_{r,p}(v_\theta)
 \frac{(h^Tv_\theta)(h_+^Tv_\theta)^f}
      {t_\theta^{(d+f)/2}}.
\end{equation}
For example, a paired-port kernel contributes
$\theta^\ell(v_\theta)_i(v_\theta)_j/t_\theta$; a port-to-loop
kernel replaces one endpoint entry by $h_+^Tv_\theta$.  The component
containing the deleted port instead has the factor $h^Tv_\theta$ at
its free endpoint.  The zero-edge free-to-loop component has these
same factors, by its full two-projector decomposition in
Section~\ref{raw:topology}.  Thus every surviving port contributes
exactly one entry of $v_\theta$, regardless of its partner.
Formula~\eqref{raw:pure-contraction} is unchanged if $v_\theta$ is
rescaled, and proves the asserted independence from pairing.

\smallskip\noindent\emph{Summation across multiplicities and parity.}
For $m>2p$, counting the two carry reflections and deleting the
ordinary edge at $p$ gives
\[
 e_C(m)=\left\lfloor\frac{m-r-1}{2}\right\rfloor+
         \left\lfloor\frac{m-p-1}{2}\right\rfloor-1.
\]
The fixed classes solve $2a\equiv m-r\pmod d$.  When $d$ is odd,
there is one such class, so $f=1$ and
$e_C=m-(r+p+5)/2$.  All factors in
\eqref{raw:pure-contraction} other than $\theta^{e_C(m)}$ are
therefore constant as $m$ increases by one.  Summing the homogeneous terms
has denominator $1-\theta$.  When $d$ is even, the values are
\[
 \begin{array}{ccc}
 &f&e_C(m)\\\hline
 m-r\text{ odd}&0&m-(r+p)/2-2\\
 m-r\text{ even}&2&m-(r+p)/2-3
 \end{array}
\]
and splitting the sum by parity gives denominator $1-\theta^2$.
These are the same homogeneous terms that occur on the step-$d$
progressions in Lemma~\ref{raw:rationality}. Formula~\eqref{raw:pure-contraction}
allows them to be regrouped across those progressions into a step-one sum
when $d$ is odd and two step-two sums when $d$ is even.
In either case the two homogeneous contributions together have a denominator
dividing $E_\tau$, where
\[
 E_j=(1-\lambda^j)(1-\mu^j)=1-S_j+c^j,\qquad
 L_d=E_d/E_\tau,\qquad \Delta=b^2-4c.
\]
Since $\tau\mid d$, the quotient $L_d$ is a polynomial: it is the
product of the two finite geometric sums for $\lambda^\tau$ and
$\mu^\tau$, which is symmetric in $\lambda,\mu$.

\smallskip\noindent\emph{Denominators from the finite path calculations.}
The recurrence has denominator $E_dR_d$, but summing the homogeneous
terms directly uses only $E_\tau$ in place of $E_d$. Could the missing
factor $L_d=E_d/E_\tau$ reappear in the coefficients of those geometric
sums? We first identify the possible denominators in a finite path
calculation, then show that none can supply a factor of $L_d$.
It is convenient to use a matrix whose entries have integer exponents.  The carry matrix is
\[
 M=\begin{pmatrix}0&z^2&z^3\\z&z^2&z^3\\z&z^2&z^3\end{pmatrix}.
\]
It is similar to $C$ and assigns the full vertex weight to the arriving
vertex of an oriented edge: $M_{ij}=z^j{\bf1}_{i+j\ge3}$.
For an oriented path with successive values $i_0,\ldots,i_\ell$, the
product
\[
 z^{i_0}\prod_{k=1}^{\ell}M_{i_{k-1},i_k}
\]
therefore gives every path vertex its full weight.  Restrict free or loop
endpoints by the appropriate indicators and sum their values.  Give every
core vertex its full weight as well, and divide the path factor by $z^i$
at each port of value $i$ shared with the core, so that each coordinate
is weighted exactly once.  For the zero-edge free-to-loop component the
same convention gives $\sum_{i=2}^3z^i$; its two-summand spectral expression
is justified by the boundary identity above.  The nonzero spectral
projector of $M$ is
\[
 \Pi_\lambda=\frac{M(M-\mu I)}{\lambda(\lambda-\mu)},
\]
with the analogous formula for $\Pi_\mu$.  Set
\[
 \mathcal A=\mathbb Q[z,(c\Delta)^{-1}],\qquad
 \mathcal S=\mathcal A[\lambda].
\]
The identities $\lambda^{-1}=\mu/c$ and
$(\lambda-\mu)^{-1}=(\lambda-\mu)/\Delta$ put both projectors in
$\operatorname{Mat}_3(\mathcal S)$.  The remaining edge and boundary
factors in the convention just specified belong to $\mathcal S$.
Weight redistribution uses integer powers of $z$, which are units in
$\mathcal A$ because $z$ divides $c\Delta$.  Constraint indicators are
$0$ or $1$; they need not be units.  Each finite amplitude is obtained
by finitely many sums and products of these factors and projector entries,
so it belongs to $\mathcal S$.  This statement concerns finite amplitudes,
not yet the sum over all multiplicities.

\smallskip\noindent\emph{Geometric denominators and conjugation.}
For the homogeneous terms, the preceding summation has geometric
denominator dividing $E_\tau$.  For a mixed assignment on any
multiplicity progression, Lemma~\ref{raw:rationality} gives denominator
$1-\lambda^{d-j}\mu^j$ with $1\le j\le d-1$, a factor of $R_d$.
Clearing these denominators and including the finite prefix therefore gives
$E_\tau R_dF_{r,p}\in\mathcal S$.

The polynomial
$\Delta=z^3(z+1)(z^2+z+4)$ is not a square in $\mathbb Q(z)$, so
\[
 \mathcal S=\mathcal A\oplus\mathcal A\lambda.
\]
The conjugation $\sigma(\lambda)=\mu=b-\lambda$ has fixed ring
$\mathcal S^\sigma=\mathcal A$: if $a_0+a_1\lambda$ is fixed, then
$a_1(\lambda-\mu)=0$, hence $a_1=0$.
The automorphism $\sigma$ of $\mathbb Q(z)(\lambda)$ fixes
$\mathbb Q(z)$ pointwise. By Lemma~\ref{raw:rationality},
$F_{r,p}\in\mathbb Q(z)$, while $E_\tau,R_d\in\mathbb Q[z]$;
hence their product is fixed by $\sigma$. Therefore
\begin{equation}\label{raw:localization}
 E_\tau R_dF_{r,p}\in\mathcal A.
\end{equation}

\smallskip\noindent\emph{Polynomial divisibility, including multiplicities.}
On the other hand, \eqref{raw:initial-data} writes
$F_{r,p}=N/(E_dR_d)$ with $N\in\mathbb Q[z]$.
Equation~\eqref{raw:localization} says precisely that $N/L_d$ belongs
to $\mathcal A$.  To conclude polynomial divisibility, we check
$\gcd(L_d,c\Delta)=1$.  At $z=0,-1$ both eigenvalues vanish and
$E_d=1$.  At the other discriminant zeros, $z^2+z+4=0$, the eigenvalues
coincide at $b/2=-2z$ and have modulus four, so again $E_d\ne0$.
Since $L_d$ divides $E_d$, it has no zero at any zero of $c\Delta$.
Write $N/L_d=G/(c\Delta)^k$ in \eqref{raw:localization}.  Then
$N(c\Delta)^k=GL_d$, and coprimality implies $L_d\mid N$ in
$\mathbb Q[z]$, with every root multiplicity.  This argument does not
require $L_d$ and $R_d$ to be coprime.  Finally
\[
 E_\tau=(1-z-z^2)K_\tau
\]
gives \eqref{raw:correction} with $A_{r,p}=N/L_d$.
Thus it is the excess factor $L_d=E_d/E_\tau$ that cancels from the
recurrence denominator.  The all-$\lambda$ and all-$\mu$ contributions
remain, and the mixed factor $R_d$ remains in the displayed denominator.

\smallskip\noindent\emph{Numerator degree.}
For the degree bound, the maximum allowed weight of a graph word is
\[
 8+4(r-1)+3(m-r-2)=3m+r-2.
\]
In the paired-factor formula for $Q_d(y,z)$, the coefficients of
$y$ and $y^2$ have degrees at most $3d$ and $4d\le6d$, respectively;
the optional middle factor has $y$-coefficient of degree $2d\le3d$.
Consequently $\deg q_j\le3dj$.  A summand in $N_s$ has degree at most
\[
 3dj+3\bigl(m_s+(k-j)d\bigr)+r-2
 \le D_0:=3d^2+3d+6p+r-2,
\]
since $m_s\le2p+d$ and $k\le d$.
Each paired factor of $Q_d(1,z)$ has degree $4d$, and the optional
middle factor has degree $2d$, giving
$\deg Q_d(1,z)=2d(d+1)$.  Thus the finite prefix in
\eqref{raw:initial-data}, after multiplication by this denominator,
contributes degree at most
\[
 2d(d+1)+6p+r-2\le D_0.
\]
Finally $\deg E_j=4j$, so cancelling $L_d$ removes degree
$4(d-\tau)$.  Substituting $r=p-d$ yields
\[
 \deg A_{r,p}\le3d^2-2d+7p+4\tau-2\le3d^2+7p+6.
\]
\end{proof}

\subsection{Uniform absolute coefficient norms}

To turn the reduced rational expression into a coefficient bound uniform
in $r,p$, we control its numerator and reciprocal denominator separately.
The correction $(1-z-z^2)F_{r,p}$ can have negative coefficients, so the
next estimates use absolute coefficient norms rather than positivity.

Set $q=(\sqrt5-1)/2$, and for a formal series $f$ write
$\|f\|_q=\sum_{j\ge0}|[z^j]f|q^j$, allowing the value infinity.
This norm is submultiplicative whenever the right side is finite.

\begin{lemma}[Uniform denominator bounds]\label{raw:wiener}
For every $d\ge1$,
\begin{equation}\label{raw:norms}
 \|B_d\|_q<10,\qquad \|B_d^{-1}\|_q<21.
\end{equation}
\end{lemma}
\begin{proof}
At $q$ the nonzero carry eigenvalues are $1,-q^2$.  Put $u=q^2$.
The polynomials $S_n$ have nonnegative coefficients, since $b,-c$ do.
For a paired mixed factor, $j\ge1$ and $n=d-2j\ge1$, its nonconstant
coefficient norm is at most
\[
 u^jS_n(q)+u^d\le u^j+u^{d-j}+u^d
                  \le a_j:=u^j+u^{j+1}+u^{2j+1}.
\]
Let
\[
 a_1=u+u^2+u^3,\qquad
 T=\sum_{j\ge2}a_j=\frac{(1+u)u^2}{1-u}+\frac{u^5}{1-u^2}.
\]
Both $a_1,T$ are less than one.  An optional central factor
$1-c^{d/2}$ costs at most $1+u$ in norm and $(1-u)^{-1}$ in
inverse norm.  Geometric-series bounds and the finite-product inequality
$\prod(1-a_j)\ge1-\sum a_j$ give
\begin{equation}\label{raw:mixed-norms}
 \|R_d\|_q\le\frac{(1+u)(1+a_1)}{1-T}<\frac{10}{3},\qquad
 \|R_d^{-1}\|_q\le\frac1{(1-u)(1-a_1)(1-T)}<6.
\end{equation}
These comparisons follow on replacing $q$ by the larger rational number
$q_+=309017/500000$, since the displayed upper bounds increase with $u$.

For completeness the fixed-factor reciprocal estimates can be certified
by a short, precisely specified finite calculation.  For either coefficient
list
\[
 (k_0,\ldots,k_e)=(1,1,1)\quad\hbox{or}\quad
 (1,1,2,1,0,-1,-1),
\]
put $v_0=1$, $v_n=-\sum_{j=1}^{\min(e,n)}k_jv_{n-j}$ for
$1\le n\le80$, $V(z)=\sum_{n=0}^{80}v_nz^n$, and $E=1-KV$.
The deliberately rounded rational inequalities in Table~\ref{tab:raw-reciprocal-residuals} are sufficient:

\begin{table}[ht]
\centering
\caption{Rational bounds for the degree-$80$ reciprocal polynomials
and their residuals at $q_+=309017/500000$.}
\label{tab:raw-reciprocal-residuals}
\begin{tabular}{ccc}
 & $\|V\|_{q_+}$ & $\|E\|_{q_+}$\\\hline
 $K_1$ & $<212/100$ & $<1/10^{16}$\\
 $K_2$ & $<343/100$ & $<1/10^{8}$
\end{tabular}
\end{table}

They are finite rational inequalities for the recurrence just displayed:
all coefficients of $E$ through degree $80$ vanish and only the next
$e$ coefficients enter its norm.  In particular $\|E\|_{q_+}<1$,
so the identity $K^{-1}=V\sum_{j\ge0}E^j$ proves, including its
infinite tail,
\[
 \|K_1^{-1}\|_q<7/2,\qquad \|K_2^{-1}\|_q<7/2.
\]
Also $\|K_1\|_q=2$ and $\|K_2\|_q=4-2q<14/5$.
Combining these with \eqref{raw:mixed-norms} proves \eqref{raw:norms}.
\end{proof}

\begin{lemma}[A finite numerator bound and signed conversion]\label{raw:signed}
Suppose $Z_{m,r,p}(q)\le B(r,p)$ for every $m>p$.  Then, with
$\mathcal N_{r,p}=3d^2+7p+6$,
\begin{equation}\label{raw:numerator-norm}
 \|H_{r,p}\|_q\le420\mathcal N_{r,p} B(r,p),
\end{equation}
and, for every $g\ge0$,
\begin{equation}\label{raw:fib-conversion}
 [z^g]F_{r,p}\le\frac94\|H_{r,p}\|_q\Fib_g.
\end{equation}
\end{lemma}
\begin{proof}
Only coefficients through degree $\mathcal N_{r,p}$ of $F_{r,p}$ can enter the polynomial
$A_{r,p}=B_d(1-z-z^2)F_{r,p}$.  Positivity of each graph polynomial
and $g\ge m+5$ imply
\[
 \sum_{j=0}^{\mathcal N_{r,p}}[z^j]F_{r,p}\,q^j\le \mathcal N_{r,p} B(r,p),
\]
since fewer than $\mathcal N_{r,p}$ multiplicities contribute.  Using
$\|1-z-z^2\|_q=2$ and \eqref{raw:norms} gives
$\|A_{r,p}\|_q\le20\mathcal N_{r,p}B(r,p)$ and then
\eqref{raw:numerator-norm}.  This step never substitutes the potentially
divergent value $F_{r,p}(q)$ for its finite truncation.

Write $H_{r,p}=\sum_jh_jz^j$.  Its support starts at degree at least
eight.  From $(1-z-z^2)^{-1}=\sum_{n\ge0}\Fib_{n+1}z^n$,
\[
 [z^g]F_{r,p}=\sum_{j=8}^g h_j\Fib_{g-j+1}.
\]
For $8\le j\le g$, Binet's formula and
$\Fib_n\le q^{1-n}$ give
\[
 \frac{\Fib_{g-j+1}}{\Fib_g}
 \le \frac{\sqrt5}{1-q^{16}}q^j<\frac94q^j.
\]
For the last strict inequality it suffices to square and check
$5<[(9/4)(1-q_+^{16})]^2$.  Take absolute values of the $h_j$.
For $g<8$ the raw coefficient is zero, which also handles $g=0$.
\end{proof}

\begin{remark}\label{raw:not-positive}
The smallest pivot pair illustrates both the Fibonacci factor and the
sign of the correction. For $r=1$, $p=2$ and $m>4$, the carry
decomposition in Section~\ref{raw:topology} gives a single path on
$m-3$ coordinates, from a free endpoint to a loop. Along it the letters
are $1,2,3$, adjacent ones are forbidden, and the final letter is $2$
or $3$. At $m=4$ this is one loop-constrained coordinate; at $m=3$
the path is empty. Thus summing over $m$ sums over every such path word,
once for each length. Read from the free end: each $1$ must be followed
by $2$ or $3$. The word splits uniquely into blocks $2,3,12,13$.
Including the two pivots gives
\[
 F_{1,2}(z)=\frac{z^8}{1-(z^2+z^3)(1+z)}
           =\frac{z^8}{(1-z-z^2)(1+z+z^2)}.
\]
In particular $(1-z-z^2)F_{1,2}$ has coefficient $-1$ at genus nine,
even though $F_{1,2}$ counts words with nonnegative coefficients.
\end{remark}

\subsection{Strong runs at the Fibonacci point}

Lemma~\ref{raw:signed} reduces the coefficient estimate to a bound on
individual graph values, uniform in the multiplicity.  We now obtain
such a bound, with decay in the pivot positions, by grouping consecutive
strong edges and estimating their interactions with carry paths.

For the remainder of the proof all matrices are evaluated at $q$.
On four coordinates let
\[
 D_{ij}=q^{(i+j)/2}{\bf1}_{i+j\ge4},\qquad
 C_{ij}=q^{(i+j)/2}{\bf1}_{i+j\ge3}{\bf1}_{i,j\le3}.
\]
Let $E$ be $D$ with its fourth row and column set to zero.  Thus $E$
is the alphabet-three strong transfer.  Put
\begin{gather}\label{raw:constants}
 R=\frac{917}{1000},\quad k_0=\frac{10003}{10000},\quad
 \chi=\frac{47597}{50000},\quad t=\frac{95761}{100000},\nonumber\\
 u_*=(1,13/50)^T,\quad v_*=(823/1250,97/625)^T,\quad
 \nu=v_*^Tu_*=\frac{10918}{15625}.
\end{gather}

\begin{lemma}[Uniform transfer block]\label{raw:block}
Let $Q=(v,w)$ be an orthonormal basis for the nonzero eigenspaces of
$C$, in the order $1,\mu=-q^2$.  For all $k,n\ge1$,
\begin{equation}\label{raw:block-bound}
 \left|Q^TD^kQ\begin{pmatrix}1&0\\0&\mu^n\end{pmatrix}\right|
 \le R^{\max(k-3,0)}u_*v_*^T
\end{equation}
entrywise.  Moreover $\|D\|_2<R$, $\|E\|_2<19/25$, and
$Q^TEQ=Q^TDQ$.
\end{lemma}
\begin{proof}
The nonsymmetric matrix similar to the cap-$a$ restriction of $D$ has
entries $q^j{\bf1}_{i+j\ge4}$.  Its products with the positive vectors
$(10000,17164,23907,23907)^T$ for $a=4$ and $(1,2,3)^T$ for $a=3$
are componentwise strictly smaller than $R$ and $19/25$ times those
vectors, respectively.  Each comparison reduces to $q^2=1-q$.
The positive-vector bound on spectral radius, and symmetry of the
original matrices, give the asserted operator norms.

One may choose
\[
 v=\frac{(q,\sqrt q,q,0)^T}{\sqrt{1+q^2}},\qquad
 w=\frac{(1,-q^{3/2},-q^2,0)^T}{\sqrt{1+q^2}}.
\]
Their required moments, writing
$\alpha_k=v^TD^kv$, $o_k=v^TD^kw$, and $\beta_k=w^TD^kw$, are
\begin{align*}
 \alpha_1=\alpha_2&=(7-6q)/5,&
 o_1&=(1-3q)/5,& o_2&=(11-18q)/5,\\
 \beta_1&=(-7+11q)/5,&
 \beta_2&=(3-4q)/5,& \beta_4&=(-181+293q)/5.
\end{align*}
These exact identities give the rational bounds
\[
 \begin{gathered}
 \alpha_1,\alpha_2<A=\frac{823}{1250},\quad
 |o_1|<\frac{171}{1000},\quad |o_2|<\frac1{40},\\
 |\beta_1|<\frac{41}{1000},\quad
 \beta_2<B=\frac{66}{625},\quad \beta_4<J=\frac{17}{1000}.
\end{gathered}
\]
Also $|\mu|<\rho=191/500$.  For $k\ge3$, Cauchy--Schwarz gives
\[
 |\alpha_k|\le R^{k-2}A,\quad
 |o_k|\le R^{k-3}\sqrt{AJ},\quad
 |\beta_3|\le\sqrt{BJ},\quad
 |\beta_k|\le R^{k-4}J\quad(k\ge4).
\]
For example, the middle bound follows from
$(Dv)^TD^{k-3}(D^2w)$.
To check all powers, write $M_*=u_*v_*^T$.  For $k=1,2$ use the
listed moments, multiplying the second column by at most $\rho$.
For $k\ge3$ the sufficient rational comparisons are
\[
 RA<M_{*,11},\quad AJ<M_{*,21}^2,\quad
 \rho^2AJ<M_{*,12}^2,\quad
 \rho^2BJ<M_{*,22}^2,\quad \rho J/R<M_{*,22}.
\]
Direct substitution proves them.  These inequalities prove
\eqref{raw:block-bound} for every power, rather than only for a finite
list of lengths.  The fourth coordinates of $v,w$ vanish, so
$Q^TEQ=Q^TDQ$.
\end{proof}

\begin{proposition}[Uniform graph-value estimate]\label{raw:pointwise}
Define
\[
 e_D=\left\lfloor\frac{r-1}{2}\right\rfloor+
     \left\lfloor\frac{p-1}{2}\right\rfloor-{\bf1}_{p\ne2r},
 \qquad
 e_E=\max\left(0,\left\lfloor\frac{p-1}{2}\right\rfloor-r\right),
\]
and $s_4(z)=z+z^2+z^3+z^4$.  For every $m>p$,
\begin{equation}\label{raw:value}
 Z_{m,r,p}(q)\le k_0^3q^8s_4(q)R^{e_D+2e_E}\chi^{d-16}.
\end{equation}
\end{proposition}
\begin{proof}
Call a maximal nonempty sequence of ordinary strong edges in a component
of the reflection graph a \emph{strong run}.  We first count its edges
and boundaries, then estimate its transfer product.

There are $\lfloor(r-1)/2\rfloor$ ordinary strong pairs targeting $r$
and $\lfloor(p-1)/2\rfloor$ targeting $p$.  The only such pair incident
to a pivot is $\{r,p-r\}$, unless $p=2r$, when this is a loop at $r$.
Deleting the pivots therefore leaves exactly $e_D$ ordinary strong edges.
Among these, a pair with both endpoints greater than $r$ must target $p$
and is uniquely specified by
\[
 r<i<p-i,\qquad
 i=r+1,\ldots,\left\lfloor\frac{p-1}{2}\right\rfloor.
\]
There are exactly $e_E$ such edges.  Neither endpoint can have a strong
edge targeting $r$, so each of these is an entire one-edge strong run.
It uses the alphabet-three matrix $E$.  Every other run of $k$ edges
is bounded by $D^k$, enlarging any smaller vertex alphabet as needed.
Every carry endpoint is greater than its target and has alphabet three;
the zero fourth row and column of $C$ enforce this alphabet at each
strong/carry junction.

For the boundary count, each nonpivot coordinate has one reflection
incidence for each target.  Such an incidence is either an ordinary
edge, a loop, or an edge lost on deleting a pivot.  Pivot $r$ has only
the possible nonloop neighbor $p-r$, and pivot $p$ only the possible
nonloop neighbor $m+r-p$; the pivots are not adjacent.  Hence there
are at most two free boundary incidences.  For each target $t$, loops
can occur only at $2i=t$ or $2i=m+t$, so there are at most four loop
incidences in total.  If $P$ is the number of path components, $F$ the
number of free incidences, and $L$ the number of loop incidences, then
\[
 2P=F+L,\qquad F\le2,\qquad L\le4,\qquad P\le3.
\]
An isolated vertex counts as a zero-edge path with two boundary
incidences, so it is included in this identity.

Let $h_4=(q^{1/2},q,q^{3/2},q^2)^T$ and let $h_{4,+}$ be the same
vector with its first entry zero.  Free endpoint vectors are bounded
by $h_4$, and loop endpoint vectors by $h_{4,+}$, because either loop
threshold forces a value at least two.  Thus all boundary norms
combined cost at most
\begin{equation}\label{raw:boundary-budget}
 \|h_4\|_2^F\|h_{4,+}\|_2^L
 =s_4(q)^{F/2}\bigl(s_4(q)-q\bigr)^{L/2}\le s_4(q),
\end{equation}
since $s_4(q)>1$ and $s_4(q)-q<1$.  The same formula handles the two
restrictions on an isolated vertex through the inner product of its
two boundary vectors.

Every ordinary cycle contains both strong and carry edges.  Indeed,
target labels alternate around a cycle.  In an all-strong cycle the
two target reflections compose to a translation by $d$ or $-d$;
in an all-carry cycle they do the same.  Successive pairs of edges
use the same ordering of target labels, so iterating that nonzero
translation cannot close the cycle.

For each cycle, rotate the transfer product to alternate nonempty
strong and carry runs.  Write $T_i=D^{k_i}$ for an ordinary strong
run and $T_i=E$ for an alphabet-three strong run, and put
$\Lambda_n=\operatorname{diag}(1,\mu^n)$.  Its partition is bounded by
\[
 \operatorname{tr}\bigl(T_1C^{n_1}\cdots T_jC^{n_j}\bigr)
 =\operatorname{tr}\prod_{i=1}^j(Q^TT_iQ\Lambda_{n_i}),
 \qquad n_i\ge1.
\]
By Lemma~\ref{raw:block}, the absolute value of each block is bounded
entrywise by $a_iM_*$, where $M_*=u_*v_*^T$ and
\[
 a_i=R^{\max(k_i-3,0)}\quad(T_i=D^{k_i}),\qquad
 a_i=1\quad(T_i=E).
\]
The latter follows from $Q^TEQ=Q^TDQ$.  Since
$M_*^j=\nu^{j-1}M_*$, taking absolute values and multiplying the
nonnegative comparison matrices gives the cycle bound
\[
 \left(\prod_i a_i\right)\operatorname{tr}(M_*^j)
 =\left(\prod_i a_i\right)\nu^j.
\]
In particular, no dimension factor is charged for a cycle.

For a path, remove the strong run at either endpoint if present and
charge its operator norm.  A path consisting entirely of strong edges
has just one such run. There are at most two endpoint runs per path.
Their operator norms account for their edges. Unless the
whole path has been removed, what remains begins and ends with carry
runs.  If it contains $j\ge1$ strong runs, its product after projection
into the carry image has the form
\[
 \Lambda_{n_0}(Q^TT_1Q\Lambda_{n_1})\cdots
                 (Q^TT_jQ\Lambda_{n_j}),\qquad n_i\ge1.
\]
Orthogonal projection and $\Lambda_{n_0}$ have operator norm at most
one.  Taking entrywise absolute values preserves the Euclidean norm
of each endpoint vector.  The same block comparison consequently bounds
this middle product by
\[
 \left(\prod_i a_i\right)\|M_*^j\|_2
 =\left(\prod_i a_i\right)\kappa\nu^j,
 \qquad \kappa=\frac{\|u_*\|_2\|v_*\|_2}{\nu}<k_0.
\]
If there is no remaining strong run, the carry product has norm at most
one.  This includes a zero-edge remainder; a path removed in its entirety
was already bounded by its endpoint norms and strong-edge charges.
Thus in every case one factor $k_0$ per path suffices, at most $k_0^3$
in total. Each strong run in a cycle or in the interior of a path
contributes the additional factor $\nu/R^3$ derived below.

To count these runs, consider the strong graph before deleting $r$.
On its vertices $1,\ldots,p-1$, the substitution $x=p-i$ turns the
reflections into $p-x$ and $p+d-x$, with ports $1,\ldots,d$.
The residue-class path argument in Section~\ref{raw:topology} applies
here too: it only requires each residue class to be nonempty, which
holds because $p-1\ge d$.  If $f_s$ classes are fixed by the residue
involution, then $f_s\le2$, because $2a\equiv p\pmod d$ has at most
two solutions.  The $(d-f_s)/2\ge d/2-1$ paired-class components
each contain an ordinary strong edge, even if their classes have only
one vertex each.  The deleted pivot $r$ is port $d$.  Removing that
endpoint can destroy at most one of these components and cannot split
another.  At least $d/2-2$ nonempty strong components therefore survive.
These components are precisely strong runs after carry edges are
restored. At most $2P\le6$ of them are endpoint runs. Thus the
number $j$ of remaining strong runs satisfies
\begin{equation}\label{raw:paid-run-count}
 j\ge d/2-8.
\end{equation}
This argument involves only the strong interval and the boundary count;
in particular it applies throughout $p<m\le2p$ as well as for $m>2p$.

The constants satisfy the exact rational comparisons
\[
 \frac{\nu}{R^3}<\chi^2<1,\qquad
 \frac{\|u_*\|_2\|v_*\|_2}{\nu}<k_0,\qquad
 \frac{19}{25}<R^3,\qquad \chi<t,\qquad R<t^2.
\]
Charge $R$ per ordinary strong edge outside the $E$ runs and $R^3$
per $E$ edge. The operator norm covers each endpoint run:
$\|D^k\|_2\le R^k$ or $\|E\|_2<19/25<R^3$.
For each remaining strong run, the preceding product estimates give
\[
 R^{\max(k-3,0)}\nu\le R^k\frac{\nu}{R^3}
 \quad\hbox{for a $D$ run},\qquad
 \nu=R^3\frac{\nu}{R^3}\quad\hbox{for an $E$ run}.
\]
Every strong edge is counted exactly once, so the total strong charge is
$R^{e_D-e_E+3e_E}=R^{e_D+2e_E}$.
By \eqref{raw:paid-run-count}, the extra factor is bounded by
\[
 (\nu/R^3)^j\le\chi^{2j}\le\chi^{d-16}.
\]
When $d<16$, the last negative exponent simply weakens the bound.
Multiplying this estimate by the path factor $k_0^3$, the global
boundary bound \eqref{raw:boundary-budget}, and the pivot weight $q^8$
proves \eqref{raw:value}.
\end{proof}

\subsection{Summation over omitted pivots}

\begin{proof}[Proof of Theorem~\ref{raw:tail}]
Combining Lemma~\ref{raw:signed} with Proposition~\ref{raw:pointwise}
gives the coefficient multiplier
\[
 945(3d^2+7p+6)k_0^3q^8s_4(q)
                 R^{e_D+2e_E}\chi^{p-r-16}.
\]
Replace $d^2$ by $p^2$.  Put $J=\lfloor(p-1)/2\rfloor$ and
$a=R^2\chi/t<1$.  Since
$R^{\lfloor(r-1)/2\rfloor}\le t^{r-2}$, the region $r\le J$
is bounded, after reversing its geometric sum, by
\[
 \sum_{r\le J}R^{e_D+2e_E}\chi^{p-r}
 \le\frac{\chi^{p-J}(Rt)^J}{Rt^2(1-R^2\chi/t)}.
\]
Here replacing the pivot-deletion indicator by one only enlarges the
bound.  For $r>J$ there are no $E$ edges; summing backwards from
$r=p-1$ gives
\[
 \sum_{J<r<p}R^{e_D}\chi^{p-r}
 \le\frac{R^Jt^p}{Rt^2(t/\chi-1)}.
\]
An empty region is harmless because the displayed upper bound is positive.

For a quadratic $P(k)=ak^2+bk+c$, define
\begin{equation}\label{raw:geometric-tail}
 \mathcal T_x(K;P)=x^K\left[
 \frac{aK^2+bK+c}{1-x}
 +\frac{(2aK+b)x}{(1-x)^2}
 +\frac{ax(1+x)}{(1-x)^3}\right]
 =\sum_{k\ge K}P(k)x^k.
\end{equation}
Set
\[
 X=Rt\chi,\quad Y=Rt^2,\quad U=1-R^2\chi/t,\quad V=t/\chi-1,
\]
and
\[
 P_o(k)=12k^2+26k+16,\qquad P_e(k)=12k^2+38k+32.
\]
These are $3p^2+7p+6$ for $p=2k+1,2k+2$, respectively.
In both cases $J=k$.  Pairing these two parities, and replacing $q$
by $q_+$ in the positive prefactor, bounds the multiplier for every
pivot $p>2K$ by the following explicit rational number:
\begin{equation}\label{raw:explicit-tail}
 \mathcal B_K=
 \frac{945q_+^8s_4(q_+)k_0^3\chi^{-16}}{Rt^2}
 \left[
 \frac{\mathcal T_X(K;\chi P_o+\chi^2P_e)}{U}
 +\frac{\mathcal T_Y(K;tP_o+t^2P_e)}{V}
 \right].
\end{equation}
All entries in this expression are the specified rational constants;
$0<X<Y<1$ and $U,V>0$.
Substitution in \eqref{raw:geometric-tail} and integer cross-multiplication
give the exact inequality
\[
 \mathcal B_{150}<\frac{221714122}{10^9}
                       <\frac{221715}{10^6}.
\]
This finite arithmetic comparison bounds the whole infinite geometric
sum in \eqref{raw:explicit-tail}.  Together with the preceding structural
lemmas it proves \eqref{raw:main} for every genus.  For a fixed genus
the pivot sum itself is finite by Lemma~\ref{raw:support}; the displayed
majorant additionally converges absolutely.
\end{proof}

\begin{corollary}[Integer envelope]\label{raw:integer-envelope}
Let $C_g(300)=\sum_{2\le p\le300}\sum_{r<p}[z^g]F_{r,p}$ be the
exact coefficient of the raw cap-four prefix.  Then the full raw
cap-four coefficient is at most
\[
 C_g(300)+{\bf1}_{g\ge604}
              \left\lfloor\frac{221715\Fib_g}{1000000}\right\rfloor.
\]
\end{corollary}
\begin{proof}
The omitted coefficient is an integer, so the floor of its real upper
bound is still an upper bound. Lemma~\ref{raw:support} makes the
omitted coefficient zero for $g<604$, which justifies the indicator.
\end{proof}

\section{Convergent masses of the genuine residual populations}
\label{mass:section}

Let $A_g$ count genuine depth-four Kunz words of genus $g$ with at
least two coordinates equal to four, and let $B_g$ count genuine depth-five
words with at least two coordinates equal to five. Write
$A(z)=\sum_g A_gz^g$ and $B(z)=\sum_g B_gz^g$.
At a positive real argument $z$, the \emph{mass} of one of these
populations is its generating-series value: each word of genus $g$
contributes the weight $z^g$. Neither population is restricted to leaves
of the semigroup tree or to words containing a coordinate equal to one.
We prove
\begin{equation}\label{mass:values}
 A(1247/2000)<2024,\qquad B(633/1000)<113.
\end{equation}
These bounds serve both remaining branches of the proof. Together they
enter the genuine-population comparison in
Proposition~\ref{mass:tail} for $g\ge837$. The bound on $B$ also supplies
$V_g$ in \eqref{cmp:V} from genus $301$ onward, within the finite
comparison of Section~\ref{sec:computation}.

Qualitative convergence already follows from Zhu's Theorem~1.5
\cite{Zhu}. Indeed, $A_g,B_g\le\widehat n_g$, where $\widehat n_g$
counts genuine words of depth at least four. Zhu bounds
$\limsup_g\widehat n_g^{1/g}$ by $r_{1.54}<1.55$. Since
$z_4=1247/2000$ and $z_5=633/1000$ are both less than $1/1.55$,
the root test gives convergence at these two arguments. This inference supplies neither
numerical constant in \eqref{mass:values}. The present section derives
the explicit masses needed for the comparison.

We use three structural facts from Section~\ref{raw:topology}: the carry-path decomposition,
its boundary conventions, and the homogeneous-pairing calculation
\eqref{raw:pure-contraction}. The numerical omitted-pivot bound is not
an input to these mass estimates. The full raw depth-four series cannot
supply the first bound in~\eqref{mass:values}: even its single
nonnegative summand $F_{1,2}$ from Remark~\ref{raw:not-positive}
diverges at $z=1/\varphi$, and hence at the larger argument $z_4$.
We therefore impose additional self-pair inequalities satisfied by every
genuine word.

\subsection{The multiple-five mass}

We first bound the depth-five mass directly by summing a geometric
majorant over both pivot positions and the multiplicity. This gives the
bound used in \eqref{cmp:V} as well as in the eventual comparison.
For $t\in\{4,5\}$ let
\[
 D_t(z)_{ij}=z^{(i+j)/2}{\bf1}_{i+j\ge t},\qquad 1\le i,j\le t,
\]
and put $s_5=\sum_{j=1}^5z^j$, $h_5=z^3+z^4+z^5$ and
$h_4=z^2+z^3+z^4$. The Schatten fourth norm of a real symmetric matrix
$M$ is $(\operatorname{tr}M^4)^{1/4}$.

\begin{lemma}\label{mass:five}
Suppose $0<a_5<b_4<1$, $\|D_5(z)\|_4\le a_5$,
$\|D_4(z)\|_4\le b_4$, $s_5\ge b_4$, $h_5\le a_5$ and
$h_4\le b_4$. Set $k=(a_5/b_4)^{1/2}$. Then
\[
 B(z)\le \frac{z^{10}(s_5/b_4)k^{-1}}
 {(1-b_4)(1-b_4k)(1-a_5)}.
\]
In particular $B(633/1000)<113$.
\end{lemma}

\begin{proof}
Fix the two rightmost fives at $r<p<m$, delete these fixed vertices,
and retain the two reflection constraints. Carry endpoints have alphabet
$\{1,2,3,4\}$; enlarging the strong-endpoint alphabet to
$\{1,2,3,4,5\}$ can only increase the partition. The ordinary graph is a
union of paths and even cycles, the latter having length at least four.
Regard an isolated vertex as a path with two endpoint incidences. Let $P$
be the number of paths, $u$ the number of free incidences, $L_5,L_4$ the
numbers of strong and carry loop incidences, and $S,C$ the numbers of
ordinary strong and carry edges. Reflection counting gives
\[
 2P=u+L_5+L_4,\quad S+C=m-3-P,\quad u\le2,
 \quad 2S+L_5\ge r+p-4.
\]
The final inequality includes $p=2r$: deleting the second pivot then
removes a strong loop rather than an ordinary strong edge.
Cauchy--Schwarz on each path bounds its endpoint factors by
$s_5^{u/2}a_5^{L_5/2}b_4^{L_4/2}$ and its edge factors by
$a_5^S b_4^C$. For a cycle of length $\ell\ge4$, Schatten H\"older and
$\|M\|_\ell\le\|M\|_4$ give the same edge bound with no dimension
factor. Alphabet projections are contractions and preserve these bounds.
Consequently the nonpivot partition is at most
\[
 b_4^{m-3}(s_5/b_4)^{u/2}k^{2S+L_5}
 \le b_4^{m-3}(s_5/b_4)k^{r+p-4}.
\]
Now put $r=1+i$, $p=2+i+j$, $m=3+i+j+\ell$, with
$i,j,\ell\ge0$, and sum the three geometric series.

At $z=633/1000$ take
\[
 a_5=\frac{876716}{10^6},\quad b_4=\frac{982645}{10^6},\quad
 \tau=\frac{928171}{10^6}.
\]
Every norm comparison is an exact rational inequality, since
\[
 \operatorname{tr}D_t(z)^4=
 \sum_{i_1,\ldots,i_4=1}^t
 z^{i_1+i_2+i_3+i_4}
 \prod_{j=1}^4{\bf1}_{i_j+i_{j+1}\ge t},\qquad i_5=i_1.
\]
The two trace sums are strictly below $a_5^4,b_4^4$, respectively;
the endpoint inequalities hold, and
$1/2<\sqrt{a_5b_4}<\tau<1$.
The preceding bound equals
$z^{10}s_5/[\sqrt{a_5b_4}(1-\sqrt{a_5b_4})(1-a_5)(1-b_4)]$.
The function $1/[x(1-x)]$ increases for $1/2<x<1$, so replacing the
square root by $\tau$ gives the rational upper bound
\[
 \frac{z^{10}s_5}{\tau(1-\tau)(1-a_5)(1-b_4)}
 <112.194409824952<113.
\]
All summations can first be made over finite boxes; monotone convergence
then proves convergence as well as the stated bound.
\end{proof}

\subsection{Conditional self-pair contraction}

For depth four we first compare two weighted populations on each fixed
graph. Write $Z_{m,r,p}(z)$ for the raw cap-four partition of
Definition~\ref{raw:definition}, with its two rightmost fours fixed at
$r<p$, including their factor $z^8$. A partition here means the sum of
$z^{\sum_j x_j}$ over the permitted coordinate assignments.
Let $Z^{\rm self}_{m,r,p}(z)$ be the subpartition avoiding
\[
 F_i=\{x_i=1,\ x_{2i}=3\},\qquad p<i<m/2.
\]
For example, at $m=12$ the word
\[
 (2,4,2,4,1,2,2,2,2,3,2)
\]
has weight $26$ and pivots $r=2$, $p=4$. The letter $x_5=1$ appears
in two retained carry pairs: $5+9=12+2$ and $5+11=12+4$.
Both partner values are $2$, so the carry inequalities hold. Every other
retained pair has value sum at least four. Yet $2x_5=2<x_{10}=3$:
this raw word violates the Kunz inequality at $(5,5)$, which targets
neither pivot. This is precisely the event $F_5$.

Every genuine word belongs to the subpartition avoiding all $F_i$, by
the noncarry Kunz inequality at $(i,i)$. Only these specified additional inequalities are
imposed on the raw population; we use the resulting subpartition as an
upper bound for the genuine population. The carry eigenvalue satisfies $\lambda(z_4)>1$
(Lemma~\ref{mass:block}), so the unrestricted carry-path estimate does
not supply decay as paths lengthen. The aim is to extract a factor less
than one for each of many self-pair restrictions, while controlling their
overlaps; these factors compensate for that growth in the later sums. For $0<z\le1$ put
\begin{gather*}
 a=z,\quad b=z^2+z^3,\quad s=a+b,\qquad
 \epsilon=\frac{abz^3}{s(ab+s^2)},\quad
 \delta=\frac{\epsilon}{(1+z)^2},\\
 \beta=1-\epsilon,\quad \gamma=1-\delta,\quad T=\gamma/\beta,
 \quad d=p-r,\quad N_{m,p}=\max\{0,\lfloor(m-1)/2\rfloor-p\}.
\end{gather*}

\begin{lemma}\label{mass:contraction}
For every pivot pair and multiplicity,
\begin{equation}\label{mass:feedback}
 Z^{\rm self}_{m,r,p}(z)
 \le \beta^{N_{m,p}} T^{1+\lfloor d/3\rfloor}Z_{m,r,p}(z).
\end{equation}
For $m\le2p$ the stronger bound $Z^{\rm self}_{m,r,p}\le Z_{m,r,p}$
will be used.
\end{lemma}

\begin{proof}
Assume $m>2p$, and let
\[
 V=\{p+1,\ldots,\lfloor(m-1)/2\rfloor\},\qquad |V|=N_{m,p}.
\]
If the raw partition is zero, so is the self-restricted partition and
there is nothing to prove. Otherwise give the raw words the probability
measure proportional to their weights. For a center $i\in V$, write
\[
 \ell_r(i)=m+r-i,\quad \ell_p(i)=m+p-i,\quad
 B_i=\{i,2i,\ell_r(i),\ell_p(i)\}.
\]
We call $B_i$ the local coordinate block, $2i$ its target, and
$\ell_r(i),\ell_p(i)$ the reflection neighbors of its center $i$.
Every coordinate in $B_i$ exceeds $p$, so its alphabet is
$\{1,2,3\}$. Both reflection neighbors exceed $m/2$.
The center has precisely these two retained neighbors, and has no loop.
All retained constraints touching the block are carry constraints of
threshold three. The only possible identification in $B_i$ is
$2i=\ell_r(i)$ or $2i=\ell_p(i)$.

\emph{Local compatibility with a raw exterior.}
We first work with the raw population, before conditioning on avoidance
of any other $F_j$. Suppose first that the four vertices are distinct.
Fix the coordinates outside
$B_i$ to an assignment admitting at least one completion. The family
\[
 x_i=1,\qquad x_{2i}=3,\qquad
 x_{\ell_r(i)},x_{\ell_p(i)}\in\{2,3\}
\]
has total block weight $ab^2z^3$. The only retained edges at $i$ join it to its two reflection
neighbors. On each, the value sum is at least $1+2=3$.
Every other retained pair touching the block has one value at least
two and the other at least one, even when the other coordinate lies
outside the block. Thus the displayed assignments extend every
feasible exterior and realize $F_i$.
Dropping all restrictions except the two center--neighbor edges bounds the
whole block partition by
\[
 s(ab^2+bs^2).
\]
Indeed setting the center to one contributes $ab^2$, allowing the
center to be two or three contributes at most $bs^2$, and the target
contributes at most $s$.
The ratio of the displayed numerator to this denominator is $\epsilon$.
If the target coincides with one reflection neighbor, the corresponding
numerator and denominator are $abz^3$ and $ab^2+bs^2$. Their ratio is
\[
 \frac{az^3}{ab+s^2}>\epsilon.
\]
These bounds also allow any additional edge or loop among the block
vertices, since the failure families satisfy it and it only reduces
the denominator. This establishes compatibility with every feasible
raw exterior. To multiply avoidance factors we must additionally check
that the chosen local assignments survive all preceding avoidances;
that is the purpose of the ordering and overlap argument below.

\emph{An ordering of the centers.}
For distinct centers draw $i\to j$ when
\begin{equation}\label{mass:dependency}
 2j\in\{\ell_r(i),\ell_p(i)\}.
\end{equation}
If $F_j$ has already been excluded and $x_j=1$, setting the shared
coordinate $x_{2j}$ to three would reintroduce $F_j$. Processing $i$
first preserves both choices at that reflection neighbor. We construct a set of postponed
centers of size at most $1+\lfloor d/3\rfloor$, and an order of all
other centers in which every arrow goes to a later or postponed center.
Put
\[
 c=(2m+r+p)/6,\qquad I=[c-d/2,c+d/2].
\]
Process centers outside $I$ in decreasing order of $|i-c|$, with
arbitrary tie breaking. For an arrow $i\to j$, write
$j=(m+t-i)/2$ with $t\in\{r,p\}$. Since $m+t-3c=\pm d/2$,
\[
 |j-c|\le |i-c|/2+d/4<|i-c|\qquad(i\notin I).
\]
Thus an outside center never points to an earlier center. The weak
version of the same inequality shows that both affine maps preserve
$I$, so no center in $I$ points to an outside center.

The remaining centers form the consecutive integer set $V\cap I$,
with at most $d+1$ elements. If there are $d+1$ elements, postpone one
endpoint. The remaining set has at most $d$ elements and diameter
strictly less than $d$. For any fixed target vertex $j$, its two
possible incoming sources $m+r-2j$ and $m+p-2j$ differ by $d$;
therefore its indegree in the remaining graph is at most one.
Every nontrivial strongly connected component is consequently a simple
directed cycle: strong connectivity supplies an internal incoming edge
at every vertex, leaving exactly one internal incoming and outgoing
edge at each vertex. A two-cycle would give
\[
 i+2j=m+t_1,\quad j+2i=m+t_2,
 \qquad i-j=t_2-t_1\in\{-d,0,d\},
\]
which is impossible for distinct vertices of this set. Self-edges were
excluded in \eqref{mass:dependency}. Every directed cycle thus has at
least three vertices. Postpone one vertex from each cyclic component
and put the remaining acyclic graph in topological order, with $i$
before $j$ on every edge $i\to j$. The postponed set has size
\[
 f\le 1+\lfloor d/3\rfloor.
\]
All centers not postponed have now been ordered; process the postponed
ones last in any order.

\emph{Survival under the preceding avoidances.}
Consider a nonpostponed center $i$, and condition on avoidance of all
preceding $F_j$. On a fixed exterior assignment outside $B_i$ with
positive mass under that conditioning, we claim that the entire local
failure family above still survives. A preceding event disjoint from
$B_i$ is already avoided by the fixed exterior assignment. For one
meeting $B_i$, the possible intersections are exhaustive as follows.
An earlier center cannot be a reflection neighbor, since centers are
below $m/2$ and reflection neighbors above it. If $j=2i$, its center
has value three, so $F_j$ is false. If $2j=i$, its target has value
one, so $F_j$ is false.
Equal centers or equal doubled targets force $j=i$. The only remaining
possibility is a preceding target at a reflection neighbor, precisely an
arrow $i\to j$, which the ordering excludes. In a three-vertex block,
a target at the shared neighbor $2i$ would again force $j=i$.

It follows that the conditional numerator on this exterior slice is
at least the full weight already displayed. The conditional denominator
is at most the unrestricted block denominator displayed there.
The exterior weight cancels. Averaging over the exterior slices proves
\[
 \Pr(F_i\mid\text{all preceding avoidances})\ge\epsilon.
\]
This argument does not require preceding avoidances to be measurable
outside $B_i$; their survival was checked on each proposed completion.

For a postponed center use the single assignment
\[
 x_i=1,\qquad x_{2i}=3,\qquad
 x_{\ell_t(i)}=2\quad\text{when }\ell_t(i)\ne2i.
\]
It also makes false every preceding event meeting the block, including
an event whose target is a remaining reflection neighbor, now of value
two.
Its weight in the four-vertex case is $z^8$, giving conditional
probability at least
\[
 \frac{z^8}{s(ab^2+bs^2)}=\delta.
\]
In the three-vertex case its weight is $z^6$, and the bound is larger
because $s>z^2$. The same slice argument therefore applies in any order
to the postponed centers. Multiplying the conditional avoidance bounds
gives
\[
 Z^{\rm self}_{m,r,p}\le\beta^{N_{m,p}-f}\gamma^f Z_{m,r,p}
 \le\beta^{N_{m,p}} T^{1+\lfloor d/3\rfloor}Z_{m,r,p}.
\]
If any conditioning event has zero mass, the final left side is zero
and the assertion follows immediately. For $m\le2p$ there are no
eligible centers, and the raw bound suffices.
\end{proof}

\subsection{A uniform envelope for the absolute spectral expansion}

We separate the all-positive mode because it has a two-step multiplicity
progression; the remaining modes gain an additional decay factor on the
longer progressions used below. These relations will allow the complete
multiplicity sums to be bounded using finitely many bases.

From now on set $z=z_4=1247/2000$, and define
\begin{gather*}
 D_{ij}=z^{(i+j)/2}{\bf1}_{i+j\ge4},\quad 1\le i,j\le4,\\
 C_{ij}=z^{(i+j)/2}{\bf1}_{i+j\ge3}{\bf1}_{i,j\le3},\quad
 J_3=\operatorname{diag}(1,1,1,0),\quad E=J_3DJ_3.
\end{gather*}
The nonzero carry eigenvalues are $\lambda>0>\mu$, where
$\lambda^2-b\lambda-ab=0$ and $\mu=b-\lambda$.
Let $Q=(q_+,q_-)$ have orthonormal eigenvectors as columns, with
$q_+$ nonnegative, and put $\rho=|\mu|/\lambda$. Thus
\[
 C^n=\lambda^nq_+q_+^T+\mu^nq_-q_-^T\qquad(n\ge1).
\]

We specify the spectral partitions used below. For a fixed reflection
graph with $m>2p$, retain each maximal strong run as an exact matrix product,
including its endpoint alphabet projections, and expand every maximal
carry run of positive length by the displayed identity. A mode
assignment chooses $+$ or $-$ for each carry run. After all strong blocks and boundaries have
been contracted, let $W_\omega$ denote its scalar weight, including
$z^8$ and any component having no carry edge. Define
\begin{equation}\label{mass:mode-definitions}
 P_{m,r,p}=W_{(+,\ldots,+)},\qquad
 A^{\rm abs}_{m,r,p}=\sum_\omega|W_\omega|,\qquad
 M_{m,r,p}=A^{\rm abs}_{m,r,p}-P_{m,r,p}.
\end{equation}
The positive mode has nonnegative endpoint and strong-block factors,
so $P_{m,r,p}\ge0$. Taking absolute values of the projected blocks and
endpoint factors before multiplying gives exactly $A^{\rm abs}_{m,r,p}$.
In particular
\begin{equation}\label{mass:spectral-majorization}
 0\le P_{m,r,p}\le A^{\rm abs}_{m,r,p},\qquad
 0\le Z_{m,r,p}\le P_{m,r,p}+M_{m,r,p}=A^{\rm abs}_{m,r,p}.
\end{equation}
Here $P_{m,r,p}$ is the all-positive spectral contribution. The
\emph{mixed} remainder $M_{m,r,p}$ is the sum of absolute weights of
all assignments containing a negative mode, including the all-negative
assignment. These are groups of spectral terms, not disjoint word
populations. This convention differs from the separation into
all-$\lambda$, all-$\mu$, and genuinely mixed terms used in
Section~\ref{raw:topology}. Below, ``pure'' refers to the all-positive
contribution only.

There is one carry component of zero ordinary-edge length that must
still receive two modes. It occurs at
$m=3p-2r>2p$, when deletion of port $p$ leaves the isolated vertex
$i=2p-r$ with a carry loop. Its free and loop boundary vectors are
\[
 h=(z^{1/2},z,z^{3/2},0)^T,\qquad
 h_+=(0,z,z^{3/2},0)^T.
\]
Both lie in $\operatorname{im}C$, so its scalar weight is split as
\begin{equation}\label{mass:zero-edge}
 h^Th_+=(h^Tq_+)(q_+^Th_+)+(h^Tq_-)(q_-^Th_+).
\end{equation}
This is a scalar identity with the specified boundary weights, rather
than an identity expressing the full matrix $C^0$ as two projectors.
It is the convention for $W_\omega$ at that corner; its absolute
mode sum is at most $\|h\|_2\|h_+\|_2$ by Cauchy--Schwarz.
An isolated vertex in the strong core is kept as its ordinary positive
scalar and is not given this artificial carry expansion.

The constants $R,S,L,u,v$ below are local to this section at $z=z_4$;
they are distinct from the bounds used in Section~\ref{raw:section} at
$q=\varphi^{-1}$. Set
\begin{equation}\label{mass:constants}
 R=\frac{23389}{25000},\quad S=\frac{72789}{100000},\quad
 L=\frac{101778011}{10^8},\quad
 u=\begin{pmatrix}1\\7/25\end{pmatrix},\quad
 v=\begin{pmatrix}47/50\\9/50\end{pmatrix}.
\end{equation}

\begin{lemma}\label{mass:block}
We have $\|D\|_2<R$, $\|E\|_2<S$, $1<\lambda<L$ and
$\rho<39/100$. For every $k,n\ge1$, entrywise,
\begin{equation}\label{mass:actual-block}
 \left|Q^TD^kQ\operatorname{diag}(\lambda^n,\mu^n)\right|
 \le R^k\lambda^nuv^T.
\end{equation}
For a single restricted strong edge and every $n\ge1$,
\[
 \left|Q^TEQ\operatorname{diag}(\lambda^n,\mu^n)\right|
 \le S\lambda^nuv^T.
\]
\end{lemma}

\begin{proof}
Conjugating $D$ by $H=\operatorname{diag}(z^{j/2})_{j=1}^4$ gives
$M:=H^{-1}DH=(z^j{\bf1}_{i+j\ge4})_{ij}$. The positive test vector
$(1000000,1710962,2377416,2377416)^T$ satisfies $Mx<Rx$
componentwise. For the corresponding three-dimensional restriction $M_3$,
$(1000000,2146345,3002947)^T$ satisfies $M_3x<Sx$. Perron--Frobenius
and symmetry prove the norm bounds. The quadratic equation gives the
stated $\lambda,\rho$ inequalities by rational squaring.

Write $\alpha_k=q_+^TD^kq_+$, $\omega_k=q_+^TD^kq_-$ and
$\nu_k=q_-^TD^kq_-$. The following rational strict upper bounds on
absolute values suffice:
\[
\begin{array}{cccc}
 k&|\alpha_k|&|\omega_k|&|\nu_k|\\\hline
 1&68/100&18/100&5/100\\
 2&69/100&3/100&11/100\\
 4&\text{not needed}&\text{not needed}&2/100
\end{array}
\]
They can be checked without numerical eigenvectors. Set
\[
\begin{aligned}
 y_+&=(b/\lambda,1,1,0)^T,& y_-&=(\lambda/a,-1,-1,0)^T,\\
 \kappa_+&=a(b/\lambda)^2+b,& \kappa_-&=\lambda^2/a+b.
\end{aligned}
\]
Then $q_\pm=Hy_\pm/\sqrt{\kappa_\pm}$. Substitute the positive root in
these finite expressions, enclose its square root by rational squares,
and compare; the denominator signs are positive. This also gives the
same first-power moments for $E$, because $J_3q_\pm=q_\pm$.

For $k\ge2$, Cauchy--Schwarz gives
$|\alpha_k|/R^k\le(69/100)/R^2$. For $k\ge3$ it gives
\[
 \frac{|\omega_k|}{R^k}\le
 \frac{\sqrt{(69/100)(2/100)}}{R^3},\qquad
 \frac{|\nu_3|}{R^3}\le
 \frac{\sqrt{(11/100)(2/100)}}{R^3},
\]
and for $k\ge4$, $|\nu_k|/R^k\le(2/100)/R^4$.
Indeed apply $D$ to $q_+$ and $D^2$ to $q_-$ for the off-diagonal
estimate, and $D^2$ to both copies of $q_-$ for the last estimate.
For $k=1,2$ use the table directly. The second column in
\eqref{mass:actual-block}, after division by $\lambda^n$, gains
$\rho^n\le39/100$. Comparison with the four entries of $uv^T$ proves
the assertion; square the positive off-diagonal bounds if desired.
The same comparisons with denominator $S$ prove the restricted case.
This extends the finite moment checks to every power.
\end{proof}

The hypothesis $n\ge1$ is essential: the assertion is about a strong
run followed by a positive-length carry run. For example,
$\nu_2/R^2>3/25>(uv^T)_{22}$, so a bare strong-run assertion would
be false. The zero-edge component in \eqref{mass:zero-edge} is handled by its
boundary estimate; \eqref{mass:actual-block} is never used with $n=0$.

Here are the endpoint and component consequences of the lemma. A vertex
with index greater than $r$ has at most one strong incidence, so every
internal vertex of a strong run has the full alphabet $\{1,2,3,4\}$.
Alphabet restrictions on such a run therefore occur only at its endpoints.
At a carry junction they do not change the projected block, since
$J_3Q=Q$. Thus a run between carry components has the projected matrix
$Q^TD^kQ$; if both endpoints of a strong edge exceed $r$, that edge is
an entire one-edge strong run and the sharper $E$ estimate applies.
At a free or loop endpoint, absorb the alphabet and unary restrictions
into the boundary vector. Strong prefixes and suffixes are then bounded
by ordinary operator norms, using $\|E\|_2<S$ for the restricted
one-edge case. Alphabet projections and $Q^T$ are contractions, while
taking coordinatewise absolute values preserves Euclidean norm. Moreover
\[
 \|u\|_2^2\|v\|_2^2=\frac{77173}{78125}<1,\qquad
 v^Tu=\frac{619}{625}<1.
\]
Every ordinary cycle has both a strong and a carry edge: two successive
reflections in a cycle of either single kind translate an integer by the
nonzero distance $p-r$, and cannot close. Thus cycles can be written as
products of blocks with $k,n\ge1$. Products introduce no extra path
factor, and cycle traces are at most one after normalization. The at
most two free endpoint incidences together contribute at most $s_4=\sum_{j=1}^4z^j$; loop endpoint
norms are at most one, since $z^2+z^3+z^4<1<s_4$.
Counting strong edges gives
\[
 e_D=\lfloor(r-1)/2\rfloor+\lfloor(p-1)/2\rfloor
       -{\bf1}_{p\ne2r},\qquad
 e_E=\max\{\lfloor(p-1)/2\rfloor-r,0\}.
\]
With $\eta=R/L$ and $\sigma=S/R$, define
\begin{equation}\label{mass:raw-envelope}
 \mathcal E_{m,r,p}=z^8s_4L^{m-3}\eta^{e_D}\sigma^{e_E}.
\end{equation}
The argument proves $Z_{m,r,p}\le\mathcal E_{m,r,p}$ for every $m>p$,
and $A^{\rm abs}_{m,r,p}\le\mathcal E_{m,r,p}$ for every $m>2p$.
For short multiplicities the ordinary graph estimate uses the same
block bounds whenever there are carry edges and ordinary operator
norms otherwise.
The number of carry edges is at most $m-3-e_D$:
the graph has $m-3$ vertices and is a union of paths and cycles. The
restricted-edge replacement supplies $(S/R)^{e_E}$, while $L>1$
permits the stated larger carry exponent.

\subsection{Finite bases and complete multiplicity tails}

Write $f_d=1+\lfloor d/3\rfloor$. Combining
\eqref{mass:feedback} with \eqref{mass:spectral-majorization} gives the
nonnegative decomposition that will be summed:
\begin{equation}\label{mass:master-sum}
 A(z)\le
 \sum_{r<p}\sum_{m=p+1}^{2p}Z_{m,r,p}
 +\sum_{r<p}\sum_{m>2p}T^{f_d}\beta^{N_{m,p}}
             (P_{m,r,p}+M_{m,r,p}).
\end{equation}
It is enough initially to read this inequality over finite ranges;
the bounds below establish convergence. We keep the exact factor
$\beta^{N_{m,p}}$ in this decomposition. For the finite bases and mixed tails
we may enlarge it using the rational normalization
\begin{equation}\label{mass:w}
 w=\frac{70651823329321825}{72057594037927936},\qquad
 w^2\ge\beta,\quad w<1,\quad J=(wL)^2<1.
\end{equation}

Multiplying every vertex weight, including the two pivots, by $w$
gives the normalized partitions
\[
 \widehat P_{m,r,p}=w^{m-1}P_{m,r,p},\quad
 \widehat A^{\rm abs}_{m,r,p}=w^{m-1}A^{\rm abs}_{m,r,p},\quad
 \widehat M_{m,r,p}=w^{m-1}M_{m,r,p}.
\]
For $m>2p$ put
\[
 K_{m,r,p}=T^{f_d}w^{-2p-((m-1)\bmod2)}.
\]
Since $m-1-2p-((m-1)\bmod2)=2N_{m,p}$, for $X=P,M$ we have
\begin{equation}\label{mass:normalized-majorization}
 T^{f_d}\beta^{N_{m,p}} X_{m,r,p}
 \le K_{m,r,p}\widehat X_{m,r,p}.
\end{equation}
A normalized block in \eqref{mass:actual-block} gains $w^{k+n}$
on both sides; the boundary vectors each gain $w^{1/2}$ and the
pivot factor gains $w^2$. In particular the normalized zero-edge
component is $w h^Th_+$, with both terms of
\eqref{mass:zero-edge} multiplied by $w$.

The carry-path decomposition of Section~\ref{raw:topology} identifies the
carry components under $m\mapsto m+2d$: each paired-port component gains
four ordinary edges and each fixed-port component gains two, including
the component from which $p$ was deleted. The strong blocks and
boundary factors remain the same. Every assignment contributing to
$M$ uses a negative mode on at least one such component. Its absolute
weight therefore gains a factor at most
$(w\lambda)^{2d}\rho^2$ after normalization. Summing these termwise
bounds proves
\begin{equation}\label{mass:mixed-progression}
 \widehat M_{m+2d,r,p}
 \le (w\lambda)^{2d}\rho^2\widehat M_{m,r,p}.
\end{equation}

For the all-positive contribution, use the homogeneous-pairing
calculation \eqref{raw:pure-contraction} with the unit vector $q_+$.
The rank-one projector $q_+q_+^T$ supplies the same vector $q_+$ at
every surviving port. Multiply these endpoint entries by the weights
and constraints of the fixed strong core, and sum over all its allowed
coordinate values, including the surviving ports. The resulting scalar
is independent of which ports were paired in the carry graph. The
carry-loop endpoints each contribute
$q_+^Th_+$, and the free endpoint left by deleting $p$ contributes
$h^Tq_+$. The number of carry loops depends only on parity, and the
total carry exponent increases by two under $m\mapsto m+2$.
The zero-edge split in \eqref{mass:zero-edge} has exactly these factors.
Consequently
\begin{equation}\label{mass:pure-progression}
 \widehat P_{m+2,r,p}=(w\lambda)^2\widehat P_{m,r,p}
 \qquad(m>2p).
\end{equation}
Thus the two bases $2p+1,2p+2$ suffice for the pure sum, whereas the
$2d$ bases $2p+s$, $1\le s\le2d$, suffice for the mixed bound.
The multiplier $K_{m,r,p}$ is constant on each progression being used.

Here are explicit definitions of the three finite pivot sums. For a
positive integer cutoff $P_0$, set
\begin{align}
 \mathsf S(P_0)
 &=\sum_{p=2}^{P_0}\sum_{r=1}^{p-1}\sum_{m=p+1}^{2p}Z_{m,r,p},
 \label{mass:short-box}\\
 \mathsf P(P_0)
 &=\sum_{p=2}^{P_0}\sum_{r=1}^{p-1}\sum_{s=1}^{2}
 \frac{K_{2p+s,r,p}\widehat P_{2p+s,r,p}}
      {1-(w\lambda)^2},\label{mass:pure-box}\\
 \mathsf M(P_0)
 &=\sum_{p=2}^{P_0}\sum_{r=1}^{p-1}\sum_{s=1}^{2d}
 \frac{K_{2p+s,r,p}\widehat M_{2p+s,r,p}}
      {1-(w\lambda)^{2d}\rho^2}.\label{mass:mixed-box}
\end{align}
All denominators are positive because $w\lambda<wL<1$.
Equations~\eqref{mass:normalized-majorization}--\eqref{mass:pure-progression}
show exactly how these finite sums bound the corresponding complete
multiplicity sums in \eqref{mass:master-sum}.

For evaluating them one may use the right eigenvectors $y_\pm$ from
the proof of Lemma~\ref{mass:block}, with left vectors
$z^j y_{\pm,j}/\kappa_\pm$. Their positive diagonal change of basis
telescopes through complete products and endpoints and commutes with
entrywise absolute values. It therefore yields the same pure and
absolute-mode scalars just defined.

It remains to sum over pivots, rather than merely checking finitely
many of them. Write
\[
 r=2a+u_0,\quad p=2b+v_0,\quad u_0,v_0\in\{1,2\},\quad
 h=\begin{cases}0,&u_0<v_0,\\1,&u_0\ge v_0.\end{cases}
\]
The exact domain is $a\ge0$, $b\ge a+h$. Also
$e_D\ge a+b-1$ and $e_E=\max\{b-2a-u_0,0\}$. Define
\[
 \mathcal H_{u_0v_0}(A,B)
 =\sum_{a\ge0}\sum_{b\ge a+h}A^aB^b
      \sigma^{\max\{b-2a-u_0,0\}}.
\]
Splitting at $b=2a+u_0$ gives
\begin{equation}\label{mass:parity-sum}
 \mathcal H_{u_0v_0}(A,B)=
 \frac{B^h}{(1-B)(1-AB)}
 -\frac{B^{u_0+1}}{(1-B)(1-AB^2)}
 +\frac{B^{u_0+1}\sigma}{(1-B\sigma)(1-AB^2)}.
\end{equation}
The first two terms are the sum of finite inner strips; they do not
require $B<1$. The convergence requirements are
\begin{equation}\label{mass:convergence}
 AB<1,\qquad AB^2<1,\qquad B\sigma<1.
\end{equation}
At $B=1$ the removable expression is read by its limit, although none
of the evaluations below has $B=1$. Let
$\mathcal H^{>P}_{u_0v_0}$ be this rational expression minus the finite
sum over $2b+v_0\le P$, $0\le a\le b-h$. This is exactly the omitted
pivot sum, not an asymptotic approximation.

The short multiplicities $p+1\le m\le2p$, summed in
\eqref{mass:raw-envelope}, contribute for $p>P$ at most
\begin{equation}\label{mass:short}
 \frac{z^8s_4L^{-2}}{\eta(L-1)}
 \sum_{u_0,v_0}\left\{
 L^{2v_0}\mathcal H^{>P}_{u_0v_0}(\eta,\eta L^4)
 -L^{v_0}\mathcal H^{>P}_{u_0v_0}(\eta,\eta L^2)\right\}.
\end{equation}
Indeed the multiplicity sum of $L^{m-3}$ is
$L^{-2}(L^{2p}-L^p)/(L-1)$, and $\eta^{e_D}\le\eta^{a+b-1}$.

For the mixed tail put
\[
 t=\frac{252061464497}{250000000000},\qquad t^3\ge T,
 \qquad c_0=(39/100)^2,
\]
and define the rational function, on its convergence domain,
\begin{equation}\label{mass:G}
 G_P(x)=\sum_{u_0,v_0}L^{2v_0}x^{v_0-u_0}
 \mathcal H^{>P}_{u_0v_0}
       (\eta/x^2,\eta L^4x^2).
\end{equation}
This is the nonnegative sum over omitted pivots of
$L^{2p}\eta^{a+b}\sigma^{e_E}x^d$.
Pair the bases $s=2j+1,2j+2$, $0\le j<d$. Their normalization leaves
$w^{s-1-((s-1)\bmod2)}$, so their envelope is
$(1+L)L^{2p-2}J^j$. Thus the mixed multiplicity factor is
\[
 H_d=\frac{1-J^d}{(1-J)(1-c_0J^d)}.
\]
The following elementary bound retains the numerator that would be
lost by replacing its geometric sum with $d$.

\begin{lemma}\label{mass:three-evaluation}
For $0<J,c_0<1$, $a_d\ge0$ and $G(x)=\sum_{d\ge1}a_dx^d<\infty$,
\[
 \sum_{d\ge1}a_dx^dH_d\le
 \frac{G(x)-G(xJ)+\dfrac{c_0}{1-c_0}
              \{G(xJ)-G(xJ^2)\}}{1-J}.
\]
\end{lemma}
\begin{proof}
For $X=J^d$, use
$1/(1-c_0X)=1+c_0X/(1-c_0X)\le1+c_0X/(1-c_0)$,
multiply by $(1-X)/(1-J)$ and sum. Both differences are sums of
nonnegative terms of convergent series.
\end{proof}

Since $T^{1+\lfloor d/3\rfloor}\le Tt^d$, the mixed omitted tail is
at most
\begin{equation}\label{mass:mixed}
 \frac{C_0}{1-J}\left[
 G_P(t)-G_P(tJ)+\frac{c_0}{1-c_0}\{G_P(tJ)-G_P(tJ^2)\}\right],
 \qquad C_0=\frac{z^8s_4(1+L)L^{-2}T}{\eta}.
\end{equation}
All parameters here are rational. For $x=t,tJ,tJ^2$ the conditions
\eqref{mass:convergence} hold: passing from $x$ to $xJ$ leaves $AB$
unchanged and decreases $AB^2$ and $B\sigma$. The initial inequalities
are checked with the exact $w$ of \eqref{mass:w}, rather than replacing
$w^2$ by the slightly smaller $\beta$.

For the pure omitted sector, return to the exact factor $\beta^{N_{m,p}}$
in \eqref{mass:master-sum}. Equations~\eqref{mass:spectral-majorization}
and \eqref{mass:raw-envelope} give, for every $m>2p$,
\[
 T^{f_d}\beta^{N_{m,p}} P_{m,r,p}
 \le T^{f_d}\beta^{N_{m,p}} A^{\rm abs}_{m,r,p}
 \le T^{f_d}\beta^{N_{m,p}}\mathcal E_{m,r,p}.
\]
This argument does not require $P_{m,r,p}\le Z_{m,r,p}$, which need
not hold for a signed spectral expansion. To bound the sum, add the
nonnegative terms $p<m\le2p$ and use
$\beta^{\lfloor(m-1)/2\rfloor-p}$ throughout. On the added range its
exponent is nonpositive; on the original range the factor equals $\beta^{N_{m,p}}$.
Put
$J_0=L^2\beta<1$,
\[
 A_0=\eta/t^2,\quad B_0=\eta J_0\beta^{-2}t^2,\qquad
 M_1=L^{-2}J_0+L^{-1},\quad M_2=(L^{-2}+L^{-1})J_0.
\]
An upper bound for this enlarged sum with $p>P$ is
\begin{equation}\label{mass:full-feedback}
 z^8s_4T\sum_{u_0,v_0}
 \frac{\beta^{-v_0}t^{v_0-u_0}M_{v_0}}{\eta(1-J_0)}
 \mathcal H^{>P}_{u_0v_0}(A_0,B_0).
\end{equation}
For $p=2b+v_0$ the exact multiplicity sum needed here is
\[
 \sum_{m\ge p+1} L^{m-3}\beta^{\lfloor(m-1)/2\rfloor-p}
 =\beta^{-v_0}(J_0\beta^{-2})^b
     \frac{M_{v_0}}{1-J_0}.
\]
Substitute the pivot bounds to obtain \eqref{mass:full-feedback}.
Again $J_0<1$ and all three conditions in \eqref{mass:convergence}
hold. The short-tail parameters in \eqref{mass:short} also satisfy
those three conditions. Thus none of these expressions subtracts
divergent infinite sums.

\subsection{The finite interval calculation}

For the finite pivot prefixes the preceding formulas require finitely
many exact graph evaluations. For each fixed graph, path dynamic
programming sums over the endpoint coordinate value, and a cycle is
evaluated by fixing its initial coordinate value and summing the closing
edge. Loop restrictions
are applied at the corresponding vertex. For pure and absolute-mode
partitions the same multiplication is performed on the projected
blocks. An upper enclosure of $\widehat A^{\rm abs}_{m,r,p}$ minus
a lower enclosure of $\widehat P_{m,r,p}$ is a valid upper enclosure
of $\widehat M_{m,r,p}$.

All interval endpoints are integers in units of $2^{-56}$. Each
multiplication and signed division is rounded outwards. The rational
$w$, a quadratic enclosure for $\lambda$, and a rational enclosure
for $T$ are checked before evaluation. Integer endpoints are guarded
by $2^{120}$; products that exceed signed 128-bit intermediate capacity
are evaluated with arbitrary-precision integers. The resulting upper
enclosures of the sums
\eqref{mass:short-box}--\eqref{mass:mixed-box} are as follows.
Here $U=2^{56}$. The input enclosures are
\[
 \frac{73338785529844359}{U}\le\lambda\le
 \frac{73338785529844360}{U},
\]\[
 \frac{73854862990399022}{U}\le T\le
 \frac{73854862990399023}{U}.
\]
The stored output fields \texttt{pure\_hi\_units},
\texttt{short\_hi\_units}, and \texttt{mixed\_hi\_units} enclose,
respectively, the sums $\mathsf P$, $\mathsf S$, and $\mathsf M$ in
\eqref{mass:pure-box}, \eqref{mass:short-box}, and
\eqref{mass:mixed-box} at the cutoffs shown below. They are integer
upper endpoints in units of $U^{-1}$; the table displays their values
after division by $U$. These are weighted sums: $\mathsf S$ covers
the short multiplicities, while the progression denominators in
$\mathsf P$ and $\mathsf M$ include their complete multiplicity tails.
They are not coefficients or counts of genuine words. The pure
calculation requires only two base graphs per pivot pair, whereas the
short and mixed calculations require a number that grows with $p$ and
$d$; Table~\ref{tab:mass-boxes} gives the graph counts and bounds. We use a larger pivot cutoff for
the pure part and bound the respective omitted ranges separately in
Table~\ref{tab:mass-assembly}.

\begin{table}[htbp]
\centering
\caption{Finite pivot sums at $z=1247/2000$ with complete multiplicity
bounds. The integer numerators are scaled by $U=2^{56}$.}
\label{tab:mass-boxes}
\begin{tabular}{p{0.40\linewidth}r p{0.36\linewidth}}
\hline
Domain and evaluated quantity & Graph cases & Certified upper bound\\
\hline
$\mathsf P(400)$, two pure parity bases
 &159600&$101235082714477743831/U$\\
$\mathsf S(100)$, short bases
 &333300&$3738038979785516612/U$\\
$\mathsf M(100)$, mixed bases
 &333300&$36617965435406954/U$\\
\hline
\end{tabular}
\end{table}

The counts follow directly from
$2\sum_{p=2}^P(p-1)=P(P-1)$ and
$\sum_{p=2}^P\sum_{r=1}^{p-1}p
=\sum_{p=2}^P\sum_{r=1}^{p-1}2(p-r)=P(P+1)(P-1)/3$
for the short and mixed bases separately. The displayed bounds already include the respective
geometric multiplicity denominators and self-pair multipliers, not just
the base-graph partitions. The finite graph evaluations establish the
displayed interval bounds;
the progression identities and inequalities above supply their extension
to every remaining multiplicity for these pivots.

Combining the boxes with \eqref{mass:short}, \eqref{mass:mixed} and
\eqref{mass:full-feedback} gives the following convenient outward
ceilings. Each finite decimal in this table denotes the indicated
rational number.

\begin{table}[htbp]
\centering
\caption{Rational upper bounds for the six terms covering the cap-four
mass. Every displayed decimal is a terminating rational value.}
\label{tab:mass-assembly}
\begin{tabular}{lr}
\hline
Contribution & Strict upper bound\\\hline
Pure prefix, $p\le400$ & $1404.918996618068$\\
Short prefix, $p\le100$ & $51.875711778803$\\
Mixed prefix, $p\le100$ & $0.508176354267$\\
Short omitted tail, $p>100$ & $116.098136035634$\\
Mixed omitted tail, $p>100$ & $446.302688776606$\\
Full self-pair tail, $p>400$ & $3.745874303725$\\\hline
\end{tabular}
\end{table}

The sum of the displayed ceilings is below $2024$. Direct addition
before rounding gives the slightly sharper certified value
\[
 A(1247/2000)<2023.449583867101<2024.
\]
Every genuine word has a unique pair of rightmost maxima and is
covered by the contracted population. In \eqref{mass:master-sum}, the
short part is covered by $\mathsf S(100)$ and \eqref{mass:short}, the
pure part by $\mathsf P(400)$ and \eqref{mass:full-feedback}, and the
mixed part by $\mathsf M(100)$ and \eqref{mass:mixed}. All six bounds
are finite. Monotone convergence applied to the nonnegative sums proves
convergence and completes the proof of \eqref{mass:values}.

\section{Compensation and the final genus comparison}
\label{mass:final}

The two masses in~\eqref{mass:values} bound the sizes of the
exceptional families through Lemma~\ref{cmb:compensation}. On the other
side of the comparison, Proposition~\ref{cmb:capacity} bounds the number
of unused targets from below. Proposition~\ref{mass:tail} shows that the
targets suffice from genus $837$ onward. At that genus the sharper estimate
\eqref{cmb:capacity-improved} is needed.

Lemma~\ref{cmb:compensation} combines the depth-four residual leaves
with the nonleaf parents of the unique-five residual leaves into disjoint
parts of the population counted by $A$. It is a positive-real weighted
comparison, not a coefficientwise comparison of equal genera. Substituting
\eqref{mass:values} in \eqref{cmb:mass-extraction} gives, for every $g$,
\begin{equation}\label{mass:residual}
 |\mathcal E_g|<2024z_4^{-g}+113z_5^{-g},\qquad
 z_4=1247/2000,\quad z_5=633/1000.
\end{equation}

\begin{proposition}\label{mass:tail}
For every source genus $g\ge837$,
$|\mathcal E_g|<L(g)$, where $L(g)$ is the integer lower bound on
unused-target capacity in Proposition~\ref{cmb:capacity}.
\end{proposition}

\begin{proof}
Put
\[
 \bar q=\frac{309017}{500000},\qquad s_* =\frac{223607}{100000},
 \qquad r_* =\frac{809}{500}.
\]
The rational inequalities $\bar q^2+\bar q>1$, $s_*^2>5$ give
$\bar q>\varphi^{-1}$ and $s_*>\sqrt5$. Binet's formula therefore implies
\[
 \frac{|\mathcal E_g|}{\Fib_g}
 < \Xi(g):=\frac{s_*\{2024(\bar q/z_4)^g+113(\bar q/z_5)^g\}}
 {1-\bar q^{2g}}.
\]
The denominator uses the smaller $1-\varphi^{-2g}$ for both parities;
replacing $\varphi^{-1}$ by $\bar q$ only enlarges the upper bound.
The sharper estimate \eqref{cmb:capacity-improved} gives, throughout
this range,
\[
 \frac{L(g)}{\Fib_g}\ge
 C(g):=1+\frac{927}{500}
       -\frac{50r_*^2}{((63/100)r_*)^g}.
\]
Exact rational evaluation at the joining genus gives
\[
 \Xi(837)<2.850575297189
 <2.853985748250<C(837).
\]
Since $\bar q<z_4,z_5$ and $(63/100)r_*>1$, the numerator terms of
$\Xi(g)$ strictly decrease, its denominator increases, and $C(g)$
increases. The comparison follows for every $g\ge837$. Failure of
this particular scalar comparison at $836$ has no implication for
monotonicity at that genus; it is covered by the finite comparison.
\end{proof}

\section{The finite comparison}\label{sec:computation}

It remains to prove $|\mathcal E_g|\le L(g)$ for $1\le g\le836$.
We apply the residual envelope of Proposition~\ref{cmb:raw-envelope}
using two integer upper bounds: $U_g$ for the complete raw cap-four
coefficient, and $V_g$ for the genuine multiple-five count $B_g$.
The former combines retained exact coefficients with the omitted-pivot
bound; the latter switches from exact raw coefficients to the genuine
mass estimate after genus $300$. Thus the calculation compares explicit
bounds, without enumerating all numerical semigroups in this range.
All coefficients in this section are integers. The two pivots are the
rightmost occurrences of the maximum coordinate; this choice prevents
multiple counting when the pivot polynomials are summed.

Let $C_g(300)$ be the coefficient of $z^g$ in the sum of the raw
depth-four polynomials with rightmost pivot $p\le300$, over all
multiplicities and left pivots. Set
\begin{equation}\label{cmp:U}
 U_g=C_g(300)+\mathbf1_{\{g\ge604\}}
       \left\lfloor\frac{221715}{10^6}\Fib_g\right\rfloor.
\end{equation}
Corollary~\ref{raw:integer-envelope}, which includes the omitted-pivot
bound of Theorem~\ref{raw:tail}, makes $U_g$ an upper bound for the
complete raw cap-four coefficient. The support inequality for cap $a$,
\begin{equation}\label{cmp:support}
 2g\ge(a-1)m+p+1,
\end{equation}
also shows that $p>300$ cannot contribute below genus $604$ when $a=4$.
The floor is legitimate because the omitted coefficient is an integer.
The displayed tail term is activated in $U_g$ at $g=604$ and in the
shifted parent term $U_{g-1}$ at $g=605$. This support threshold does
not assert that every later omitted coefficient is nonzero.

For depth five, write $C^{(5)}_g$ for the coefficient of the complete raw
two-pivot population and define
\begin{equation}\label{cmp:V}
 V_g=\begin{cases}
 C^{(5)}_g,&0\le g\le300,\\
 \left\lfloor113(1000/633)^g\right\rfloor,&g\ge301.
 \end{cases}
\end{equation}
The stored exact depth-five array ends at degree $300$, which explains
the switch in this definition. The first line is an exact raw coefficient
calculation, hence an upper bound on the genuine count $B_g$. The second line bounds $B_g$ by
Lemma~\ref{mass:five}: positivity gives $B_gz_5^g\le B(z_5)<113$,
and integrality permits the floor. It is not a weighted-mass claim for
the larger raw population. Coefficients with negative indices are zero.

Let $N_g$ denote the no-one coefficients in \eqref{cmb:N}. With the
choices \eqref{cmp:U} and \eqref{cmp:V}, the residual bound
\eqref{cmb:Eupper} becomes
\begin{equation}\label{cmp:residual}
 E^+_g=U_g-N_g-N_{g-3}-N_{g-4}
       +U_{g-1}-N_{g-1}+V_g.
\end{equation}
The shifted parent term has no sole-final-one subtraction: such a parent is
permitted by the unique-five map.

\begin{lemma}[Exact finite comparison]\label{cmp:finite}
For every $1\le g\le836$, the integers defined above satisfy
$0\le E^+_g\le L(g)$.
\end{lemma}

\begin{proof}[Verification and coverage]
For a truncation degree $H$, inequality~\eqref{cmp:support} bounds every
contributing multiplicity. The retained depth-four calculation has
$H=897$, $p\le300$, and
\[
 3\le m\le597,\qquad
 2\le p\le\min(300,m-1),\qquad 1\le r<p.
\]
The stored array goes through degree $897$; the comparison in this
lemma uses only degrees through $836$.
These ranges comprise $17\,775\,550$ graph triples in $595$ disjoint
multiplicity classes. Each component polynomial is computed by the path or
cycle recurrence specified above, and the resulting products are summed and
truncated to give $C_g(300)$.

For depth five with $H=300$, the support inequality gives $m\le149$.
The sufficient complete run over $3\le m\le149$ and all pivot pairs visits
\[
 \sum_{m=3}^{149}\binom{m-1}{2}=\binom{149}{3}=540\,274
\]
graph triples. The resulting exact array contains all $301$ coefficients
through degree $300$. Thus no cap-five omitted-pivot estimate is used in the
exact range of~\eqref{cmp:V}.

The rational series defining $N_g$ and $L(g)$ are expanded by exact
coefficient recurrences. Substitution into
\eqref{cmp:U}--\eqref{cmp:residual} checks all $836$ inequalities. Put
$\Delta_g=L(g)-E^+_g$. The smallest integer slack is $1$, at $g=1,2$.
Among $1\le g\le836$, the smallest Fibonacci-normalized slack is
\[
 \min_g\frac{\Delta_g}{\Fib_g}
 =\frac{15028705}{33489287}>0.448761,
\]
attained at $g=42$. Table~\ref{tab:finite-transitions} separates the
$300/301$ change in the depth-five estimate from the $604/605$
activation of the depth-four tail in the two shifted terms. The
$836/837$ change between the finite and eventual arguments is checked
after the lemma. The displayed decimal entries are downward truncations of
exact rational values. The capacity-normalized slack also satisfies
\[
 \frac{\Delta_g}{L(g)}>0.179861\qquad(1\le g\le836),
\]
with its minimum attained at $g=605$.

\begin{table}[ht]
\centering\small
\caption{Transition checks for the finite comparison.}
\label{tab:finite-transitions}
\begin{tabular}{@{}rllr@{}}
\toprule
$g$ & depth-five term & depth-four tail status & $\Delta_g/\Fib_g$\\
\midrule
300 & exact $C^{(5)}_{300}$ & absent & $>0.871878$\\
301 & analytic mass bound & absent & $>0.683638$\\
603 & analytic mass bound & absent & $>0.872080$\\
604 & analytic mass bound & enters $U_g$ & $>0.650368$\\
605 & analytic mass bound & enters $U_{g-1}$ & $>0.513344$\\
836 & analytic mass bound & both terms present & $>0.513474$\\
\bottomrule
\end{tabular}
\end{table}

For completeness, the polynomial calculation introduces no numerical
approximation. A polynomial truncated at degree $H$ is evaluated at the
digit base $b=2^{H+2}$. The coefficient at actual weight $h\le H$ is at most
the number of positive compositions of $h$, hence less than $b$. Components
are placed in a fixed coordinate order, and the unique pivot choice preserves
this interpretation on summation. Minimum-weight normalization changes only
the position of a digit; a state is pruned only when its least possible
completed weight exceeds $H$. A negative normalization shift is made only
after an exact divisibility check. Discarded higher digits therefore cannot
alter a retained coefficient.

The exact transition rows, input identities, and regeneration procedure are
recorded in the repository~\cite{LiuCompanion}. Proposition~\ref{cmb:raw-envelope}
and the preceding analytic bounds identify these integer comparisons
with the required bound on $|\mathcal E_g|$.
\end{proof}

At the other end of the finite range,
$\Delta_{836}/\Fib_{836}>0.513474$. The scalar analytic comparison of
Proposition~\ref{mass:tail} fails at $836$, but at $837$ its exact rational
margin satisfies
\[
 C(837)-\Xi(837)>0.0034104510,
 \qquad \frac{C(837)-\Xi(837)}{C(837)}>0.0011949783.
\]
Thus the finite and analytic arguments meet without an uncovered genus.

\begin{proof}[Proof of Theorem~\ref{thm:main}]
For $g=0$, the two counts are $n_0=n_1=1$.
For $1\le g\le836$, the residual bound and
Lemma~\ref{cmp:finite} give $|\mathcal E_g|\le L(g)$.
For $g\ge837$, Proposition~\ref{mass:tail} gives the same comparison.
Applying~\eqref{cmb:fundamental} yields the result.
\end{proof}

\section{Concluding remarks}

The stronger Fibonacci inequality remains open:
\[
 n_{g+1}\ge n_g+n_{g-1}\qquad(g\ge1).
\]
The present construction guarantees enough unused targets
to accommodate its exceptional sources, but does not provide the
additional $n_{g-1}$ targets required for this inequality. The finer
error terms studied by Zhu~\cite{Zhu} offer a complementary way to
formulate that obstruction.

It would also be useful to replace the comparison of cardinalities
for the exceptional sources by a direct operation on their Kunz words,
with an explicit inverse on its image. The theorem already guarantees
an injection between the two finite genus classes. The question is
whether the remaining matching can be explained by a structural rule,
rather than by estimates and finite verification.

\appendix
\setcounter{table}{0}
\renewcommand{\thetable}{\thesection.\arabic{table}}
\renewcommand{\theHtable}{\thesection.\arabic{table}}
\setcounter{figure}{0}
\renewcommand{\thefigure}{\thesection.\arabic{figure}}
\renewcommand{\theHfigure}{\thesection.\arabic{figure}}
\section{Computational material}\label{app:computation}

The computational part has three distinct tasks: generating the raw
coefficient arrays, bounding the finite weighted graph sums, and checking
the resulting scalar inequalities. Table~\ref{tab:computations}
records their parameter ranges. The verification of stored outputs and
the regeneration of those outputs are separate procedures.

\begin{table}[ht]
\centering\small
\caption{Finite inputs to the proof. The retained depth-four arrays
cover a wider range than the theorem requires.}
\label{tab:computations}
\begin{tabular}{@{}p{.37\textwidth}p{.36\textwidth}r@{}}
\toprule
Calculation & Range & Graphs\\
\midrule
Raw depth-four coefficients & $p\le300$, $3\le m\le597$;
degrees $0$--$897$ & $17\,775\,550$\\
Raw depth-five coefficients & $3\le m\le149$; all pivots;
degrees $0$--$300$ & $540\,274$\\
Pure-mode weighted sum & $p\le400$ & $159\,600$\\
Short and mixed weighted sums & $p\le100$ & $666\,600$\\
\bottomrule
\end{tabular}
\end{table}

The depth-five row describes sufficient regeneration coverage.
The retained file was originally produced over $3\le m\le293$
($4\,149\,466$ parameter triples); the support inequality makes all
additional contributions zero through degree $300$. Regeneration over
$3\le m\le149$ gives the same $301$ coefficients.

The coefficient calculation uses arbitrary-precision integer arithmetic.
The weighted graph sums use rational enclosures with outward rounding,
and the final scalar comparisons use exact rational arithmetic.
The calculation of an enclosure, its propagation through a recurrence,
and the final strict comparison are all specified by the mathematical
bounds in Sections~\ref{raw:section}, \ref{mass:section}, and~\ref{sec:computation}.

\paragraph{\textit{Code and data availability.}}
The verification code, exact finite outputs and regeneration instructions
are available in the public repository~\cite{LiuCompanion}. Its verification guide maps
the computations in the text to their inputs and commands, with separate
instructions for checking stored results and regenerating coefficient
arrays and interval bounds. All Python checks use the standard library;
regeneration requires a C++17 compiler with 128-bit integer support and
the GMP library. For the finite comparison, the capacity bound is the sum
of $L_{\rm nonbin}(g)$ from the safe and early-one prefix families and
$L_{\rm bin}(g)$ from the maximum in~\eqref{cmb:binary-capacity}. The
corresponding stored fields are
\texttt{free\_safe\_lower}, \texttt{free\_early\_one\_lower} and
\texttt{free\_binary\_lower}. The historical field
\texttt{total\_free\_lower} is not this sum and is not used in the
final comparison.

\paragraph{\textit{Formal verification.}}
A Lean~4 formalization of Theorem~\ref{thm:main} is available in the
companion repository~\cite{LiuLean}. It represents actual numerical
semigroups and their genus, and formalizes the partial injection and the
finite and analytic comparisons. The formal argument uses Lean's kernel
and standard axioms, together with compiled evaluation for 279 finite
computations covering the depth-four and depth-five coefficient arrays
and the finite prefix bounds. Fresh proof-term replay retains this
compiled-evaluation trust; it does not independently re-prove those
computations by kernel reduction. The cited asymptotic result of Zhai is
background and is not a dependency of the formalized main theorem.

\section*{AI statement}
This work was developed through continuous collaboration between the
author and AI models. GPT-6 Astra (OpenAI) and Fable 5.1 (Anthropic)
contributed to proof development and checking, exact computation,
literature research, and preparation and review of the manuscript,
with further editorial review for this version by Codex and GPT-6 Pro
(OpenAI), Claude Fable 5.2 (Anthropic), and Gemini 4 Carbon (Google).
The author takes full responsibility for all mathematical results
and the contents of this manuscript.

\section*{Acknowledgments}
The author thanks Nathan Kaplan for a helpful discussion and detailed
comments on an earlier version, especially his suggestions for clarifying
the proof structure, the terminology, and the relation to previous work.

\end{document}